# A VECTOR SMALL-GAIN THEOREM FOR GENERAL NONLINEAR CONTROL SYSTEMS[§]

**Iasson Karafyllis**[*] **and Zhong-Ping Jiang**[**]

[*]**Department of Environmental Engineering,**
Technical University of Crete,
73100, Chania, Greece
email: ikarafyl@enveng.tuc.gr

[**]**Department of Electrical and Computer Engineering,**
Polytechnic Institute of New York University,
Six Metrotech Center, Brooklyn, NY 11201, U.S.A.
email: zjiang@control.poly.edu

**Abstract**
A new Small-Gain Theorem is presented for general nonlinear control systems. The novelty of this research work is that vector Lyapunov functions and functionals are utilized to derive various input-to-output stability and input-to-state stability results. It is shown that the proposed approach recovers several recent results as special instances and is extendible to several important classes of control systems such as large-scale complex systems, nonlinear sampled-data systems and nonlinear time-delay systems. An application to a biochemical circuit model illustrates the generality and power of the proposed vector small-gain theorem.

## 1. Introduction

The small-gain theorem has been widely recognized as an important tool for robustness analysis and robust controller design within the control systems community. For instance, classical small-gain theorems [4,43] have played a crucial role for linear robust control of uncertain systems subject to dynamic uncertainties [44]. As introduced in the framework of classical small-gain, an essential condition for input-output stability of a feedback system is that the loop gain is less than one. This condition relying upon on the concept of linear finite-gain was first relaxed by Hill [7] and then Mareels and Hill [26] using the notions of monotone gain and nonlinear operators. Some nonlinear versions of the classical small-gain theorem are derived for input-output stability of nonlinear feedback systems. Quickly after the birth of the notion of input-to-state stability (ISS) originally introduced by Sontag [33], a nonlinear, generalized small-gain theorem was developed in [9]. This nonlinear ISS small-gain theorem differs from classical small-gain theorems and the nonlinear small-gain theorem of [7] and [26] in several aspects. One of them is that both internal and external stability properties are discussed in a single framework, while only input-output stability is addressed in previous small-gain theorems. As demonstrated in [9] and the subsequent work of many others, nonlinear small-gain has led to new solutions to several challenging problems in robust nonlinear control, such as stabilization by partial-state and output feedback, robust adaptive tracking, and nonlinear observers. More interestingly, this perspective of nonlinear small-gain can find useful applications in monotone systems, an important class of systems in mathematical biology (see [1,2,5]). Further extensions of this tool to the cases of non-uniform in time stability, discrete-time systems and Lyapunov characterizations are pursued by several authors independently; see, for instance, [6,8,10,11,13,16,17,20,38,40].

---

[§] This work has been supported in part by the NSF grant DMS-0504462.

This paper takes a step further to broaden the applicability and generality of nonlinear small-gain by removing two essential restrictions in previous small-gain theorems. A common feature of the earlier nonlinear small-gain theorems is that the semi-group property is required implicitly or explicitly for the solutions of the feedback system in question, whether the feedback system is described by ordinary differential equations or takes the form of hybrid and switched systems. We will adopt a weak semi-group property which is much more relaxed than the semi-group property (see [18,19]). As shown in our recent work [20], the weak semi-group property allows studying a wide class of nonlinear feedback systems such as hybrid and switched systems. While the small-gain theorem to be introduced in this paper retains the above-mentioned features of our earlier nonlinear ISS small-gain theorem, it will enjoy an additional new feature, that is, both finite-dimensional and infinite-dimensional systems can be addressed. Here we mainly focus on time-delay systems described by retarded functional differential equations. We will achieve this somehow ambitious task by making use of vector Lyapunov functions and functionals. Because of the latter, we coin our new small-gain theorem presented in this paper "vector small gain". The advantage of vector Lyapunov function versus single Lyapunov function in nonlinear stability analysis has been well documented in past literature [25,28]. Recent work in [21, 24] provides further evidence on the usefulness of vector Lyapunov functions to the case of input-to-state stability. In this paper, we will show that the vector small-gain theorem can recover several newly introduced small-gain theorems for large-scale complex systems [3,15,21,42]. In addition, it is shown that, within vector small-gain, uniform and non-uniform input-to-output and input-to-state stability properties can be studied for various important classes of nonlinear dynamical controlled systems. Examples of these systems include those represented by Ordinary Differential Equations (ODEs), Retarded Functional Differential Equations (RFDEs), and sampled-data systems. An interesting application to a biochemical control circuit model is given to demonstrate the effectiveness of this vector small-gain methodology.

The rest of the paper is organized as follows. In Section 2 we provide certain useful results on monotone discrete-time systems. The results contained in this section are used extensively in subsequent sections. Section 3 of the paper provides a brief review of the system-theoretic framework introduced in [18,19,20] and the main result is stated (Theorem 3.1). In Section 4 sufficient Lyapunov-like conditions for the verification of the hypotheses of Theorem 3.1 are presented for three types of systems: (i) Systems described by ODEs, (ii) Systems described by RFDEs and (iii) Sampled-Data systems. The results contained in Section 4 are exploited in Section 5, where examples and applications of the vector small-gain methodology are given. The conclusions of the paper are provided in Section 6. The proofs of most of the results of the paper are given in the Appendix.

**Notations** Throughout this paper we adopt the following notations:

- We denote by $K^+$ the class of positive, continuous functions defined on $\Re_+ := \{x \in \Re : x \geq 0\}$. We say that a function $\rho : \Re_+ \to \Re_+$ is positive definite if $\rho(0) = 0$ and $\rho(s) > 0$ for all $s > 0$. By $K$ we denote the set of positive definite, increasing and continuous functions. We say that a positive definite, increasing and continuous function $\rho : \Re_+ \to \Re_+$ is of class $K_\infty$ if $\lim_{s \to +\infty} \rho(s) = +\infty$. By $KL$ we denote the set of all continuous functions $\sigma = \sigma(s,t) : \Re_+ \times \Re_+ \to \Re_+$ with the properties: (i) for each $t \geq 0$ the mapping $\sigma(\cdot,t)$ is of class $K$; (ii) for each $s \geq 0$, the mapping $\sigma(s,\cdot)$ is non-increasing with $\lim_{t \to +\infty} \sigma(s,t) = 0$.

- $\Re_+^n := (\Re_+)^n = \{(x_1,...,x_n)' \in \Re^n : x_1 \geq 0,...,x_n \geq 0\}$. $\{e_i\}_{i=1}^n$ denotes the standard basis of $\Re^n$. $Z_+$ denotes the set of non-negative integers.

- Let $x, y \in \Re^n$. We say that $x \leq y$ if and only if $(y - x) \in \Re_+^n$. We say that a function $\rho : \Re_+^n \to \Re^+$ is of class $N_n$, if $\rho$ is continuous with $\rho(0) = 0$ and such that $\rho(x) \leq \rho(y)$ for all $x, y \in \Re_+^n$ with $x \leq y$.

- For $t \geq t_0 \geq 0$ let $[t_0,t] \ni \tau \to V(\tau) = (V_1(\tau),...,V_n(\tau))' \in \Re^n$ be a bounded map. We define $[V]_{[t_0,t]} := \left(\sup_{\tau \in [t_0,t]} V_1(\tau),..., \sup_{\tau \in [t_0,t]} V_n(\tau)\right)$. For a measurable and essentially bounded function $x : [a,b] \to \Re^n$ $\operatorname*{ess\,sup}_{t \in [a,b]} |x(t)|$ denotes the essential supremum of $|x(\cdot)|$.

- We say that $\Gamma : \Re_+^n \to \Re_+^m$ is non-decreasing if $\Gamma(x) \leq \Gamma(y)$ for all $x, y \in \Re_+^n$ with $x \leq y$. For an integer $k \geq 1$, we define $\Gamma^{(k)}(x) = \underbrace{\Gamma \circ \Gamma \circ ... \circ \Gamma}_{k \text{ times}}(x)$, when $m = n$.

- We define $\mathbf{1} = (1,1,...,1)' \in \Re^n$. If $u \leq v$ then $\mathbf{1}u \leq \mathbf{1}v$.

- By $\| \|_\mathcal{X}$, we denote the norm of the normed linear space $\mathcal{X}$. By $| \,|$ we denote the Euclidean norm of $\Re^n$. Let $U \subseteq \mathcal{X}$ with $0 \in U$. By $B_U[0,r] := \{ u \in U \,;\, \|u\|_\mathcal{X} \leq r \}$ we denote the intersection of $U \subseteq \mathcal{X}$ with the closed sphere of radius $r \geq 0$, centered at $0 \in U$. If $U \subseteq \Re^n$ then $\mathrm{int}(U)$ denotes the interior of the set $U \subseteq \Re^n$.

- Let $U$ be a subset of a normed linear space $\mathcal{U}$, with $0 \in U$. By $\mathcal{M}(U)$ we denote the set of all locally bounded functions $u: \Re_+ \to U$. By $u_0$ we denote the identically zero input, i.e., the input that satisfies $u_0(t) = 0 \in U$ for all $t \geq 0$. If $U \subseteq \Re^n$ then $M_U$ denotes the space of measurable, locally bounded functions $u: \Re_+ \to U$.

## 2. Global Asymptotic Stability for Monotone Discrete-Time Systems

Consider the discrete-time system

$$x_{k+1} = \Gamma(x_k) \quad , \quad x_k \in \Re_+^n \tag{2.1}$$

where $\Gamma: \Re_+^n \to \Re_+^n$ is a non-decreasing map with $\Gamma(0) = 0$. For the study of the above system we adopt the standard stability notions for discrete-time systems (see for instance [12,14] and references therein). More specifically, we say that $0 \in \Re_+^n$ is a Globally Asymptotically Stable (GAS) equilibrium point for (2.1) if $\lim_{k \to \infty} \Gamma^{(k)}(x) = 0$ for all $x \in \Re_+^n$ and for every $\varepsilon > 0$ there exists $\delta > 0$ such that $|x| \leq \delta$, $x \in \Re_+^n$ implies $|\Gamma^{(k)}(x)| \leq \varepsilon$ for all $k \geq 1$. Next a necessary condition for the Global Asymptotic Stability property and a technical result that guarantees convergence to zero are provided.

**Proposition 2.1:** *If $0 \in \Re^n$ is GAS, then the following implication holds:*

$$\Gamma(x) \geq x \Rightarrow x = 0 \tag{2.2}$$

**Proof:** Suppose that there exists $x \in \Re_+^n$, $x \neq 0$ with $\Gamma(x) \geq x$. By induction, it follows that $\Gamma^{(k)}(x) \geq x$ for all $k \geq 1$. Letting $k \to \infty$ leads to a contradiction. Thus (2.2) holds. ◁

**Lemma 2.2:** *Let $\Gamma: \Re_+^n \to \Re_+^n$ be a continuous, non-decreasing map satisfying (2.2) with $\Gamma(0) = 0$. If the inequality $\Gamma(x) \leq x$ holds for some $x \in \Re_+^n$ then $\lim_{k \to \infty} \Gamma^{(k)}(y) = 0$ for all $y \in \Re_+^n$ with $y \leq x$.*

**Proof:** The sequence $\Gamma^{(k)}(x)$ is non-increasing (in the sense that $0 \leq ... \leq \Gamma^{(3)}(x) \leq \Gamma^{(2)}(x) \leq \Gamma(x) \leq x$) and therefore the limit $\lim_{k \to \infty} \Gamma^{(k)}(x) = r$ exists. By virtue of continuity of the mapping $\Gamma: \Re_+^n \to \Re_+^n$ we have $\Gamma(r) = r$. Condition (2.2) guarantees that $r = 0$ and consequently $\lim_{k \to \infty} \Gamma^{(k)}(x) = 0$. Moreover, if $y \leq x$ then by induction it holds that $\Gamma^{(k)}(x) \geq \Gamma^{(k)}(y) \geq 0$ for all $k \geq 1$, and consequently $\lim_{k \to \infty} \Gamma^{(k)}(y) = 0$. ◁

**Definition 2.3:** *Let $x = (x_1,...,x_n)' \in \Re^n$, $y = (y_1,...,y_n)' \in \Re^n$. We define $z = MAX\{x, y\}$, where $z = (z_1,...,z_n)' \in \Re^n$ satisfies $z_i = \max\{x_i, y_i\}$ for $i = 1,...,n$. Similarly for $u_1,...,u_m \in \Re^n$ we have $z = MAX\{u_1,...,u_m\}$ with $z_i = \max\{u_{1i},...,u_{mi}\}$, $i = 1,...,n$.*

**Remark 2.4:**

**(a)** The $MAX$ operator is different from the $\max$ operator (which is implicitly defined on every set as the maximal element of a set). For example, if $x = (1,2)' \in \Re^2$, $y = (0,3)' \in \Re^2$ then the operation $\max\{x, y\}$ is not defined while $MAX\{x, y\} = (1,3)'$. Of course, if $x \leq y$ then $y = \max\{x, y\} = MAX\{x, y\}$.

**(b)** If $z \leq MAX\{x, y\}$ and $y \leq MAX\{w, v\}$ then $z \leq MAX\{x, w, v\}$. Also $MAX\{x, y, w\} = MAX\{x, MAX\{y, w\}\}$. If $x \leq z$ and $y \leq z$ then $MAX\{x, y\} \leq z$.

**(c)** Let $x, y \in \Re_+^n$. In general we have $\Gamma(MAX\{x, y\}) \geq MAX\{\Gamma(x), \Gamma(y)\}$ for any non-decreasing map $\Gamma : \Re_+^n \to \Re_+^n$.

**Definition 2.5:** *We say that* $\Gamma : \Re_+^n \to \Re_+^n$ *is* **MAX-preserving** *if* $\Gamma : \Re_+^n \to \Re_+^n$ *is non-decreasing and for every* $x, y \in \Re_+^n$ *the following equality holds:*

$$\Gamma(MAX\{x, y\}) = MAX\{\Gamma(x), \Gamma(y)\} \tag{2.3}$$

The above defined MAX-preserving maps enjoy the following important property.

**Proposition 2.6:** $\Gamma : \Re_+^n \to \Re_+^n$ *with* $\Gamma(x) = (\Gamma_1(x), ..., \Gamma_n(x))'$ *is MAX-preserving if and only if there exist non-decreasing functions* $\gamma_{i,j} : \Re_+ \to \Re_+$, $i, j = 1, ..., n$ *with* $\Gamma_i(x) = \max_{j=1,...,n} \gamma_{i,j}(x_j)$ *for all* $x \in \Re_+^n$, $i = 1, ..., n$.

**Proof:** Define $\gamma_{i,j}(s) := \Gamma_i(se_j)$ for all $s \geq 0$. Let $x \in \Re_+^n$, i.e., $x = x_1 e_1 + ... + x_n e_n$ with $x_i \geq 0$, $i = 1, ..., n$. Notice that $x = MAX\{x_1 e_1, ..., x_n e_n\}$ and consequently $\Gamma(x) = MAX\{\Gamma(x_1 e_1), ..., \Gamma(x_n e_n)\}$. Therefore $\Gamma_i(x) = \max\{\Gamma_i(x_1 e_1), ..., \Gamma_i(x_n e_n)\} = \max_{j=1,...,n} \gamma_{i,j}(x_j)$. The converse statement is a direct consequence of the definition $\Gamma_i(x) = \max_{j=1,...,n} \gamma_{i,j}(x_j)$.  ◁

Next, necessary and sufficient conditions are provided for GAS of (2.1) for the case of a continuous MAX-preserving map. The proof of the following proposition is provided in the Appendix.

**Proposition 2.7:** *Suppose that* $\Gamma : \Re_+^n \to \Re_+^n$ *with* $\Gamma(x) = (\Gamma_1(x), ..., \Gamma_n(x))'$ *is MAX-preserving and there exist functions* $\gamma_{i,j} \in N_1$, $i, j = 1, ..., n$ *with* $\Gamma_i(x) = \max_{j=1,...,n} \gamma_{i,j}(x_j)$, $i = 1, ..., n$. *The following statements are equivalent:*

**(i)** $0 \in \Re^n$ *is GAS for (2.1).*

**(ii)** *It holds that* $\gamma_{i,i}(s) < s$, *for all* $s > 0$, $i = 1, ..., n$. *Furthermore, if* $n > 1$ *then the following set of small-gain conditions holds for each* $r = 2, ..., n$:

$$(\gamma_{i_1, i_2} \circ \gamma_{i_2, i_3} \circ ... \circ \gamma_{i_r, i_1})(s) < s, \ \forall s > 0$$

*for all* $i_j \in \{1, ..., n\}$, $i_j \neq i_k$ *if* $j \neq k$.

**(iii)** *The following implication holds:* $\Gamma(x) \geq x \Rightarrow x = 0$.

**(iv)** *(iii) holds and for each* $k \geq 1$ *and* $x \in \Re_+^n$ *it holds that* $\Gamma^{(k)}(x) \leq Q(x) = MAX\{x, \Gamma(x), \Gamma^{(2)}(x), ..., \Gamma^{(n-1)}(x)\}$.

**Remark 2.8:** Notice that $Q : \Re_+^n \to \Re_+^n$ is a continuous, MAX-preserving map with $Q(0) = 0$ and $Q(a) \geq a$ for all $a \in \Re_+^n$. Moreover $\Gamma(Q(x)) \leq Q(x)$ for all $x \in \Re_+^n$.

The next proposition is a useful technical result, which will be used in the following section.

**Proposition 2.9:** *Suppose that* $\Gamma : \Re_+^n \to \Re_+^n$ *with* $\Gamma(x) = (\Gamma_1(x), ..., \Gamma_n(x))'$ *is MAX-preserving and there exist functions* $\gamma_{i,j} \in N_1$, $i, j = 1, ..., n$ *with* $\Gamma_i(x) = \max_{j=1,...,n} \gamma_{i,j}(x_j)$, $i = 1, ..., n$. *Moreover, suppose that implication (2.2) holds and that* $x \leq MAX\{a, \Gamma(x)\}$ *for certain* $x, a \in \Re_+^n$. *Then* $x \leq Q(a)$, *where* $Q(a) = MAX\{a, \Gamma(a), \Gamma^{(2)}(a), ..., \Gamma^{(n-1)}(a)\}$.

**Proof:** Suppose that $x \le MAX\{a, \Gamma(x)\}$. Then $\Gamma(x) \le MAX\{\Gamma(a), \Gamma^{(2)}(x)\}$ and $x \le MAX\{a, \Gamma(a), \Gamma^{(2)}(x)\}$. By an induction argument $x \le MAX\{a, \Gamma(a), ..., \Gamma^{(k)}(a), \Gamma^{(k+1)}(x)\}$ for all $k \ge 1$. It follows from statement (iv) of Proposition 2.7 that $x \le MAX\{Q(a), \Gamma^{(k+1)}(x)\}$ for all $k \ge 1$. Since $\lim_{k \to \infty} \Gamma^{(k)}(x) = 0$, we obtain $x \le Q(a)$. ◁

## 3. A Vector Small-Gain Theorem for a Wide Class of Systems

### 3.A. Review of the System-Theoretic Framework

In this work we make use of the system theoretic framework presented in [18,19,20]. For reasons of completeness the basic notions are recalled here.

**The notion of a Control System-Definition 2.1 in [20]:** *A control system $\Sigma := (\mathcal{X}, \mathcal{Y}, M_U, M_D, \phi, \pi, H)$ with outputs consists of*

(i) *a set $U$ (control set) which is a subset of a normed linear space $\mathcal{U}$ with $0 \in U$ and a set $M_U \subseteq \mathcal{M}(U)$ (allowable control inputs) which contains at least the identically zero input $u_0$,*
(ii) *a set $D$ (disturbance set) and a set $M_D \subseteq \mathcal{M}(D)$, which is called the "set of allowable disturbances",*
(iii) *a pair of normed linear spaces $\mathcal{X}, \mathcal{Y}$ called the "state space" and the "output space", respectively,*
(iv) *a continuous map $H : \Re_+ \times \mathcal{X} \times U \to \mathcal{Y}$ that maps bounded sets of $\Re_+ \times \mathcal{X} \times U$ into bounded sets of $\mathcal{Y}$, called the "output map",*
(v) *a set-valued map $\Re_+ \times \mathcal{X} \times M_U \times M_D \ni (t_0, x_0, u, d) \to \pi(t_0, x_0, u, d) \subseteq [t_0, +\infty)$, with $t_0 \in \pi(t_0, x_0, u, d)$ for all $(t_0, x_0, u, d) \in \Re_+ \times \mathcal{X} \times M_U \times M_D$, called the set of "sampling times"*
(vi) *and the map $\phi : A_\phi \to \mathcal{X}$ where $A_\phi \subseteq \Re_+ \times \Re_+ \times \mathcal{X} \times M_U \times M_D$, called the "transition map", which has the following properties:*

1) **Existence:** *For each $(t_0, x_0, u, d) \in \Re_+ \times \mathcal{X} \times M_U \times M_D$, there exists $t > t_0$ such that $[t_0, t] \times (t_0, x_0, u, d) \subseteq A_\phi$.*

2) **Identity Property:** *For each $(t_0, x_0, u, d) \in \Re_+ \times \mathcal{X} \times M_U \times M_D$, it holds that $\phi(t_0, t_0, x_0, u, d) = x_0$.*

3) **Causality:** *For each $(t, t_0, x_0, u, d) \in A_\phi$ with $t > t_0$ and for each $(\tilde{u}, \tilde{d}) \in M_U \times M_D$ with $(\tilde{u}(\tau), \tilde{d}(\tau)) = (u(\tau), d(\tau))$ for all $\tau \in [t_0, t]$, it holds that $(t, t_0, x_0, \tilde{u}, \tilde{d}) \in A_\phi$ with $\phi(t, t_0, x_0, u, d) = \phi(t, t_0, x_0, \tilde{u}, \tilde{d})$.*

4) **Weak Semigroup Property:** *There exists a constant $r > 0$, such that for each $t \ge t_0$ with $(t, t_0, x_0, u, d) \in A_\phi$:*

(a) *$(\tau, t_0, x_0, u, d) \in A_\phi$ for all $\tau \in [t_0, t]$,*
(b) *$\phi(t, \tau, \phi(\tau, t_0, x_0, u, d), u, d) = \phi(t, t_0, x_0, u, d)$ for all $\tau \in [t_0, t] \cap \pi(t_0, x_0, u, d)$,*
(c) *if $(t + r, t_0, x_0, u, d) \in A_\phi$, then it holds that $\pi(t_0, x_0, u, d) \cap [t, t + r] \ne \emptyset$.*
(d) *for all $\tau \in \pi(t_0, x_0, u, d)$ with $(\tau, t_0, x_0, u, d) \in A_\phi$ we have $\pi(\tau, \phi(\tau, t_0, x_0, u, d), u, d) = \pi(t_0, x_0, u, d) \cap [\tau, +\infty)$.*

**The BIC and RFC properties-Definition 2.2 in [20]:** *Consider a control system $\Sigma := (\mathcal{X}, \mathcal{Y}, M_U, M_D, \phi, \pi, H)$ with outputs. We say that system $\Sigma$*

(i) *has the "Boundedness-Implies-Continuation" (BIC) property if for each $(t_0, x_0, u, d) \in \Re_+ \times \mathcal{X} \times M_U \times M_D$, there exists a maximal existence time, i.e., there exists* $t_{\max} := t_{\max}(t_0, x_0, u, d) \in (t_0, +\infty]$, such that $A_\phi = \bigcup_{(t_0, x_0, u, d) \in \mathfrak{R}_+ \times \mathcal{X} \times M_U \times M_D} [t_0, t_{\max}) \times \{(t_0, x_0, u, d)\}$. In addition, if $t_{\max} < +\infty$ then for every $M > 0$ there exists $t \in [t_0, t_{\max})$ with $\|\phi(t, t_0, x_0, u, d)\|_\mathcal{X} > M$.

**(ii)** is **robustly forward complete (RFC) from the input** $u \in M_U$ if it has the BIC property and for every $r \geq 0$, $T \geq 0$, it holds that

$$\sup\{\|\phi(t_0 + s, t_0, x_0, u, d)\|_\mathcal{X} ; u \in \mathcal{M}(B_U[0, r]) \cap M_U, s \in [0, T], \|x_0\|_\mathcal{X} \leq r, t_0 \in [0, T], d \in M_D\} < +\infty$$

**The notion of a robust equilibrium point-Definition 2.3 in [20]:** Consider a control system $\Sigma := (\mathcal{X}, \mathcal{Y}, M_U, M_D, \phi, \pi, H)$ and suppose that $H(t, 0, 0) = 0$ for all $t \geq 0$. We say that $0 \in \mathcal{X}$ is a **robust equilibrium point from the input** $u \in M_U$ for $\Sigma$ if

**(i)** for every $(t, t_0, d) \in \mathfrak{R}_+ \times \mathfrak{R}_+ \times M_D$ with $t \geq t_0$ it holds that $\phi(t, t_0, 0, u_0, d) = 0$.

**(ii)** for every $\varepsilon > 0$, $T, h \in \mathfrak{R}_+$ there exists $\delta := \delta(\varepsilon, T, h) > 0$ such that for all $(t_0, x, u) \in [0, T] \times \mathcal{X} \times M_U$, $\tau \in [t_0, t_0 + h]$ with $\|x\|_\mathcal{X} + \sup_{t \geq 0}\|u(t)\|_\mathcal{U} < \delta$ it holds that $(\tau, t_0, x, u, d) \in A_\phi$ for all $d \in M_D$ and

$$\sup\{\|\phi(\tau, t_0, x, u, d)\|_\mathcal{X} ; d \in M_D, \tau \in [t_0, t_0 + h], t_0 \in [0, T]\} < \varepsilon$$

Next we present the Input-to-Output Stability property for the class of systems described previously (see also [9,37]).

**The notions of IOS, UIOS, ISS and UISS-Definition 2.5 in [20]:** Consider a control system $\Sigma := (\mathcal{X}, \mathcal{Y}, M_U, M_D, \phi, \pi, H)$ with outputs and the BIC property and for which $0 \in \mathcal{X}$ is a robust equilibrium point from the input $u \in M_U$. Suppose that $\Sigma$ is RFC from the input $u \in M_U$. If there exist functions $\sigma \in KL$, $\beta \in K^+$, $\gamma \in N$ such that the following estimate holds for all $u \in M_U$, $(t_0, x_0, d) \in \mathfrak{R}_+ \times \mathcal{X} \times M_D$ and $t \geq t_0$:

$$\|H(t, \phi(t, t_0, x_0, u, d), u(t))\|_Y \leq \sigma(\beta(t_0)\|x_0\|_\mathcal{X}, t - t_0) + \sup_{t_0 \leq \tau \leq t} \gamma(\|u(\tau)\|_\mathcal{U}) \quad (*)$$

then we say that $\Sigma$ satisfies the **Input-to-Output Stability (IOS) property from the input** $u \in M_U$ with gain $\gamma \in N$. Moreover, if $\beta(t) \equiv 1$ then we say that $\Sigma$ satisfies the **Uniform Input-to-Output Stability (UIOS) property from the input** $u \in M_U$ with gain $\gamma \in N$.

*For the special case of the identity output mapping, i.e., $H(t, x, u) := x$, the (Uniform) Input-to-Output Stability property from the input $u \in M_U$ is called (Uniform) Input-to-State Stability ((U) ISS) property from the input $u \in M_U$.*

The reader should notice that other equivalent definitions of the ISS property are available in the literature (see [6,31]).

### 3.B. A New Small-Gain Theorem

We consider an abstract control system $\Sigma := (\mathcal{X}, \mathcal{Y}, M_U, M_D, \phi, \pi, H)$ with the BIC property for which $0 \in \mathcal{X}$ is a robust equilibrium point from the input $u \in M_U$. Suppose that there exist maps $V_i : \mathfrak{R}_+ \times \mathcal{X} \times U \to \mathfrak{R}^+$, with $V_i(t, 0, 0) = 0$ for all $t \geq 0$ ($i = 1, ..., n$) and a MAX-preserving continuous map $\Gamma : \mathfrak{R}_+^n \to \mathfrak{R}_+^n$ with $\Gamma(0) = 0$ such that the following hypotheses hold:

**(H1)** There exist functions $\sigma \in KL$, $\nu, c \in K^+$, $\zeta, a, p^u \in N_1$, $p \in N_n$, $L : \mathfrak{R}_+ \times \mathcal{X} \to \mathfrak{R}^+$ with $L(t, 0) = 0$ for all $t \geq 0$, such that for every $(t_0, x_0, u, d) \in \mathfrak{R}_+ \times \mathcal{X} \times M_U \times M_D$ the mappings $t \to V(t) = \left(V_1\left(t, \phi(t, t_0, x_0, u, d), u(t)\right), \ldots, V_n\left(t, \phi(t, t_0, x_0, u, d), u(t)\right)\right)'$ and $t \to L(t) = L\left(t, \phi(t, t_0, x_0, u, d)\right)$ are locally bounded on $[t_0, t_{\max})$ and the following estimates hold for all $t \in [t_0, t_{\max})$:

$$V(t) \leq MAX\left\{ \mathbf{1}\sigma\left(L(t_0), t - t_0\right), \Gamma\left(\left[V\right]_{[t_0, t]}\right), \mathbf{1}\zeta\left(\left[\|u(\tau)\|_{\mathcal{U}}\right]_{[t_0, t]}\right)\right\} \quad (3.1)$$

$$L(t) \leq \max\left\{ \nu(t - t_0), c(t_0), a\left(\|x_0\|_{\mathcal{X}}\right), p\left(\left[V\right]_{[t_0, t]}\right), p^u\left(\left[\|u(\tau)\|_{\mathcal{U}}\right]_{[t_0, t]}\right)\right\} \quad (3.2)$$

where $t_{\max}$ is the maximal existence time of the transition map of $\Sigma$.

**(H2)** $0 \in \Re^n$ is GAS for the discrete-time system (2.1).

**(H3)** There exist functions $b \in N_1$, $g \in N_n$, $\mu, \beta, \kappa \in K^+$ such that the following inequalities hold for all $(t, x, u) \in \Re_+ \times \mathcal{X} \times U$:

$$\mu(t)\|x\|_{\mathcal{X}} \leq b\left(L(t, x) + g(V(t, x, u)) + \kappa(t)\right) \text{ and } L(t, x) \leq b\left(\beta(t)\|x\|_{\mathcal{X}}\right) \quad (3.3)$$

where $V(t, x, u) = \left(V_1(t, x, u), \ldots, V_n(t, x, u)\right)'$.

For future reference, $V(t, x, u) = \left(V_1(t, x, u), \ldots, V_n(t, x, u)\right)'$ is called the vector Lyapunov function for the system $\Sigma := (\mathcal{X}, \mathcal{Y}, M_U, M_D, \phi, \pi, H)$.

The main result of the present work is stated next.

**Theorem 3.1:** *Consider system $\Sigma := (\mathcal{X}, \mathcal{Y}, M_U, M_D, \phi, \pi, H)$ with the BIC property for which $0 \in \mathcal{X}$ is a robust equilibrium point from the input $u \in M_U$ and suppose that there exist maps $V_i : \Re_+ \times \mathcal{X} \times U \to \Re^+$, with $V_i(t, 0, 0) = 0$ for all $t \geq 0$ ($i = 1, \ldots, n$) and a MAX-preserving continuous map $\Gamma : \Re_+^n \to \Re_+^n$ with $\Gamma(0) = 0$ such that hypotheses (H1-3) hold.*

*Then there exist functions $\tilde{\sigma} \in KL$ and $\tilde{\beta} \in K^+$ such that for every $(t_0, x_0, u, d) \in \Re_+ \times \mathcal{X} \times M_U \times M_D$ the following estimate holds for all $t \geq t_0$:*

$$V(t) \leq MAX\left\{ \mathbf{1}\tilde{\sigma}\left(\tilde{\beta}(t_0)\|x_0\|_{\mathcal{X}}, t - t_0\right), G\left(\left[\|u(\tau)\|_{\mathcal{U}}\right]_{[t_0, t]}\right)\right\} \quad (3.4)$$

*where*

$$G(s) = MAX\left\{Q\left(\mathbf{1}\sigma\left(p^u(s), 0\right)\right), Q\left(\mathbf{1}\sigma\left(p(Q(\mathbf{1}\zeta(s))), 0\right)\right), Q\left(\mathbf{1}\zeta(s)\right)\right\} \quad (3.5)$$

*and $Q(x) = MAX\left\{x, \Gamma(x), \Gamma^{(2)}(x), \ldots, \Gamma^{(n-1)}(x)\right\}$. Moreover, if $\beta, c \in K^+$ are bounded then $\tilde{\beta} \in K^+$ is bounded. Finally, if in addition to (H1-3) the following hypothesis holds:*

**(H4)** *There exists $q \in N_n$ such that the following inequality holds for all $(t, x, u) \in \Re_+ \times \mathcal{X} \times U$:*

$$\|H(t, x, u)\|_Y \leq q(V(t, x, u)) \quad (3.6)$$

*where $V(t, x, u) = \left(V_1(t, x, u), \ldots, V_n(t, x, u)\right)'$.*

then system $\Sigma$ satisfies the IOS property from the input $u \in M_U$ with gain $\gamma(s) := q\left(1 \max_{i=1,...,n} G_i(s)\right)$, where $G(s) = (G_1(s),...,G_n(s))'$ is the mapping defined by (3.5). Moreover, if $\beta, c \in K^+$ are bounded then system $\Sigma$ satisfies the UIOS property from the input $u \in M_U$ with gain $\gamma(s) := q\left(1 \max_{i=1,...,n} G_i(s)\right)$.

**Remark 3.2:** Notice that for the control input-free case, i.e., $u \equiv 0$, Theorem 3.1 implies (Uniform) Robust Global Asymptotic Output Stability (RGAOS) for the corresponding system. Moreover, if there exists $M \geq 1$ such that $\sigma(s,0) = Ms$ for all $s \geq 0$ then the functions $G_i \in N_1$ ($i=1,...,n$) with $G(s) = (G_1(s),...,G_n(s))'$ are given by:

$$G_i(s) := \varphi_i\left(\max\left\{Mp^u(s), Mp(\varphi_1(\zeta(s)),...,\varphi_n(\zeta(s))), \zeta(s)\right\}\right), \quad i=1,...,n$$

where

$$\varphi_i(s) = \max\left\{s, \max_{k=1,...,n-1} \max\left\{\left(\gamma_{i,j_1} \circ \gamma_{j_1,j_2} ... \circ \gamma_{j_{k-1},j_k}\right)(s) : (j_1,...,j_k) \in \{1,...,n\}^k\right\}\right\}, \quad i=1,...,n$$

and $\gamma_{i,j} \in N_1$, $i, j = 1,...,n$ are the functions with $\Gamma_i(x) = \max_{j=1,...,n} \gamma_{i,j}(x_j)$, $i=1,...,n$ and $\Gamma(x) = (\Gamma_1(x),...,\Gamma_n(x))'$. The reader should notice that if hypothesis (ii) of Proposition 2.7 holds for the functions $\gamma_{i,j} \in N_1$, $i, j = 1,...,n$, then for each $i=1,...,n$ we have either $\varphi_i(s) = s$ or there exists an index set $(i, j_1,..., j_k) \in \{1,...,n\}^{k+1}$ with no repeated index such that $\varphi_i(s) = \left(\gamma_{i,j_1} \circ \gamma_{j_1,j_2} ... \circ \gamma_{j_{k-1},j_k}\right)(s)$.

## 4. Vector Lyapunov Functions and Functionals

In this section we provide sufficient Lyapunov-like conditions for the verification of Theorem 3.1 for three types of systems: (i) Systems described by Ordinary Differential Equations (ODEs), (ii) Systems described by Retarded Functional Differential Equations (RFDEs) and (iii) Sampled-Data systems. Notice that since families of Lyapunov functions (or functionals) are employed, the obtained results constitute conditions for vector Lyapunov functions (or functionals) for the (U)IOS property.

### 4.A. Systems of ODEs

We consider systems described by Ordinary Differential Equations (ODEs) of the form:

$$\dot{x} = f(t,x,u,d) \quad, \quad Y = H(t,x)$$
$$x \in \Re^n, Y \in \Re^N, u \in U, d \in D, t \geq 0 \tag{4.1}$$

where $D \subseteq \Re^l$, $U \subseteq \Re^m$ with $0 \in U$ and $f: \Re_+ \times \Re^n \times U \times D \to \Re^n$, $H: \Re_+ \times \Re^n \to \Re^N$ are continuous mappings with $H(t,0) = 0$, $f(t,0,0,d) = 0$ for all $(t,d) \in \Re_+ \times D$ that satisfy the following hypotheses:

**(A1)** There exists a symmetric positive definite matrix $P \in \Re^{n \times n}$ such that for every bounded $I \subseteq \Re_+$ and for every bounded $S \subset \Re^n \times U$, there exists a constant $L \geq 0$ satisfying the following inequality:

$$(x-y)' P(f(t,x,u,d) - f(t,y,u,d)) \leq L|x-y|^2$$
$$\forall t \in I, \forall (x,u,y,u) \in S \times S, \forall d \in D$$

**(A2)** There exist $a \in K_\infty$, $\gamma \in K^+$ such that $|f(t,x,u,d)| \leq \gamma(t)a(|x|+|u|)$ for all $(t,x,u,d) \in \Re_+ \times \Re^n \times U \times D$.

**(A3)** There exist functions $V_i \in C^1(\Re_+ \times \Re^n; \Re_+)$ ($i=1,...,k$), $W \in C^1(\Re_+ \times \Re^n; \Re_+)$, $a_1, a_2, a_3, a_4 \in K_\infty$ $\mu, \beta, \kappa \in K^+$, $\zeta \in N_1$, $g \in N_k$, $\gamma_{i,j} \in N_1$, $p_i \in N_1$, $i, j = 1,...,k$, with $\gamma_{i,i}(s) \equiv 0$ for $i=1,...,k$, a family of positive definite functions $\rho_i \in C^0(\Re_+;\Re_+)$ ($i=1,...,k$) and a constant $\lambda \in (0,1)$ such that the following inequalities hold for all $(t,x,u) \in \Re_+ \times \Re^n \times U$ :

$$a_1(|H(t,x)|) \leq \max_{i=1,...,k} V_i(t,x) \leq a_2(\beta(t)|x|) \tag{4.2}$$

$$a_3(\mu(t)|x|) - g(V_1(t,x),...,V_k(t,x)) - \kappa(t) \leq W(t,x) \leq a_4(\beta(t)|x|) \tag{4.3}$$

$$\sup\left\{\frac{\partial W}{\partial t}(t,x) + \frac{\partial W}{\partial x}(t,x)f(t,d,x,u) : d \in D\right\} \leq -W(t,x) + \lambda \max\left\{\zeta(|u|), \max_{j=1,...,k} p_j(V_j(t,x))\right\} \tag{4.4}$$

and for every $i=1,...,k$ and $(t,x,u) \in \Re_+ \times \Re^n \times U$ the following implication holds:

"If $\max\left\{\zeta(|u|), \max_{j=1,...,k} \gamma_{i,j}(V_j(t,x))\right\} \leq V_i(t,x)$ then $\frac{\partial V_i}{\partial t}(t,x) + \sup_{d \in D}\frac{\partial V_i}{\partial x}(t,x)f(t,d,x,u) \leq -\rho_i(V_i(t,x))$" (4.5)

Our main result concerning systems of the form (4.1) is the following result which provides sufficient conditions for Theorem 3.1 to hold. Its proof is provided in the Appendix.

**Theorem 4.1 (Vector Lyapunov Function Characterization of the IOS property):** *Consider system (4.1) under hypotheses (A1-3). If the following set of small-gain conditions holds for each $r=2...,k$ :*

$$\left(\gamma_{i_1,i_2} \circ \gamma_{i_2,i_3} \circ ... \circ \gamma_{i_r,i_1}\right)(s) < s, \quad \forall s > 0 \tag{4.6}$$

*for all $i_j \in \{1,...,k\}$, $i_j \neq i_l$ if $j \neq l$, then system (4.1) satisfies the IOS property with gain $\gamma = a_1^{-1} \circ \theta \in N_1$ from the input $u \in M_U$, where*

$$\theta(s) := \max_{i=1,...,k} \varphi_i\left(\max\left\{\max_{i=1,...,k} p_i(\zeta(s)), \max_{i=1,...,k}\max_{j=1,...,k} \gamma_{i,j}(\varphi_j(\zeta(s))), \max_{i=1,...,k}\max_{j=1,...,k} p_i(\gamma_{i,j}(\varphi_j(\zeta(s)))), \zeta(s)\right\}\right) \tag{4.7}$$

*and*

$$\varphi_i(s) := \max\left\{s, \max_{l=1,...,k-1}\max\left\{\left(\gamma_{i,j_1} \circ \gamma_{j_1,j_2}...\circ \gamma_{j_{l-1},j_l}\right)(s) : (j_1,...,j_l) \in \{1,...,n\}^l\right\}\right\}, \quad i=1,...,k \tag{4.8}$$

*Moreover, if $\beta \in K^+$ is bounded then system (4.1) with output $Y = H(t,x)$ satisfies the UIOS property with gain $\gamma = a_1^{-1} \circ \theta \in N_1$ from the input $u \in M_U$.*

**Comment on Theorem 4.1:** The proof of Theorem 4.1 shows that inequalities (4.2), (4.5) are used for the derivation of inequalities (3.1) and (3.6), while inequalities (4.3), (4.4) are used for the derivation of inequalities (3.2) and (3.3). Hypotheses (A1) and (A2) are minimal regularity hypotheses that guarantee uniqueness of solutions and continuity of the solutions with respect to initial data for system (4.1). Moreover, it is clear from implication (4.5) that in this case the diagonal gain functions $\gamma_{i,i}(s)$ play no role in implication (4.5); that is why they are assumed to be zero.

For the ISS case where $H(t,x) = x$, one can set $W(t,x) \equiv 0$ in Theorem 4.1 and obtain the following corollary.

**Corollary 4.2 (Vector Lyapunov Function Characterization of the ISS property):** *Consider system (4.1) under hypotheses (A1-2) and suppose that there exists a family of functions $V_i \in C^1(\Re_+ \times \Re^n;\Re_+)$ ($i=1,...,k$), functions $a_1, a_2 \in K_\infty$ $\beta \in K^+$, $\zeta \in N_1$, $\gamma_{i,j} \in N_1$, $i,j=1,...,k$, with $\gamma_{i,i}(s) \equiv 0$ for $i=1,...,k$ and a family of positive definite functions $\rho_i \in C^0(\Re_+;\Re_+)$ ($i=1,...,k$), such that:*

$$a_1(|x|) \leq \max_{i=1,...,k} V_i(t,x) \leq a_2(\beta(t)|x|), \quad \forall (t,x) \in \Re_+ \times \Re^n \tag{4.9}$$

*and implication (4.5) holds for every $i=1,...,k$ and $(t,x,u) \in \Re_+ \times \Re^n \times U$. If, additionally, the small-gain conditions (4.6) hold for each $r=2...,k$ and for all $i_j \in \{1,...,k\}$, $i_j \neq i_l$ if $j \neq l$, then system (4.1) satisfies the ISS* *property with gain* $\gamma = a_1^{-1} \circ \theta \in N_1$ *from the input* $u \in M_U$, *where* $\theta \in N_1$ *is defined by (4.7), (4.8). Moreover, if* $\beta \in K^+$ *is bounded then system (4.1) satisfies the UISS property with gain* $\gamma = a_1^{-1} \circ \theta \in N_1$ *from the input* $u \in M_U$, *where* $\theta \in N_1$ *is defined by (4.7), (4.8).*

**Comparison of Theorem 4.1 and Corollary 4.2 with existing results:** The reader should compare the result of Corollary 4.2 with Theorem 3.4 in [21]. It is clear that Theorem 3.4 in [21] is a special case of Corollary 4.2 with $\gamma_{i,j}(s) = a(s)$ for all $i,j = 1,...,k$, where $a \in N_1$ with $a(s) < s$ for $s > 0$. On the other hand vector Lyapunov characterizations based on the main result in [3] (e.g., Theorem 3.6 in [21]) require Lipschitz regularity of the right-hand sides of the ODEs, since the main result in [3] is based on certain qualitative characterizations of the ISS property provided in [34,35]. This is exactly the reason why the main result in [3] does not provide a formula for the gain function of the overall system (in contrast with the main results in the present work). The same comments hold for Theorem 2 in [15]. In particular, it should be mentioned that the small-gain conditions in Corollary 4.2 are exactly the same as the cyclic small-gain conditions in [42, 15]. In order to demonstrate the applicability of our results to large-scale interconnected systems, consider the case

$$\dot{x}_i = f_i(d, x, u) \quad , \quad i = 1,...,k$$
$$x = (x_1,...,x_k)' \in \Re^N, d \in D, u \in U$$

where $x_i \in \Re^{n_i}$, $i = 1,...,k$, $N = n_1 + ... + n_k$, $D \subset \Re^l$ a non-empty compact set, $U \subseteq \Re^m$ a nonempty set with $0 \in U$, $f_i : D \times \Re^N \times U \to \Re$, $i = 1,...,k$ are locally Lipschitz mappings with $f_i(d,0,0) = 0$ for all $d \in D$, $i = 1,...,k$. We assume that the UISS property holds for each subsystem $\dot{x}_i = f_i(d, x, u)$ with input $(u, x_1,..., x_{i-1}, x_{i+1},..., x_k)$ ($i = 1,...,k$). Let $V_i \in C^1(\Re_+ \times \Re^{n_i}; \Re_+)$ ($i = 1,...,k$) be ISS-Lyapunov functions for each one of the subsystems, i.e. positive definite and radially unbounded functions for which the following inequalities hold for $i = 1,...,k$:

$$\sup_{d \in D} \left\{ \nabla V_i(x_i) f_i(d, x, u) : (d, u) \in D \times U, x = (x_1,...,x_k)' \in \Re^N, \max\left\{ \zeta(|u|), \max_{j \neq i} \gamma_{i,j}(V_j(x_j)) \right\} \leq V_i(x_i) \right\} < 0, \quad \forall x_i \neq 0$$

for certain functions $\zeta \in N_1$, $\gamma_{i,j} \in N_1$, $i,j = 1,...,k$, with $\gamma_{i,i}(s) \equiv 0$ for $i = 1,...,k$. Working with the Lyapunov functions $V_i \in C^1(\Re_+ \times \Re^{n_i}; \Re_+)$ ($i = 1,...,k$) and exploiting Corollary 4.2, we can guarantee that the UISS property holds for the above system if the small-gain conditions (4.6) hold for each $r = 2...,k$ and for all $i_j \in \{1,...,k\}$, $i_j \neq i_l$ if $j \neq l$. It should be clear that the functions $\gamma_{i,j} \in N_1$, $i,j = 1,...,k$ are the actual gain functions, i.e., the following inequalities hold for all $i = 1,...,k$, $t \geq 0$, $x(0) \in \Re^N$ and $u \in M_U$:

$$V_i(x_i(t)) \leq \max\left\{ \sigma_i(V_i(x_i(0)), t), \zeta\left(\sup_{0 \leq \tau \leq t}|u(\tau)|\right), \max_{j \neq i} \gamma_{i,j}\left(\sup_{0 \leq \tau \leq t} V_j(x_j(\tau))\right) \right\}$$

for certain $\sigma_i \in KL$ ($i = 1,...,k$), which are nothing else but the inequalities of the max-formulation of the UISS property for each subsystem $\dot{x}_i = f_i(d, x, u)$ with input $(u, x_1,..., x_{i-1}, x_{i+1},..., x_k)$ ($i = 1,...,k$). Again it can be seen that in this case the diagonal gain functions $\gamma_{i,i}(s)$ play no role; that is why they are assumed to be zero.

**4.B. Systems Described by RFDEs**

Let $D \subseteq \Re^l$ be a non-empty set, $U \subseteq \Re^m$ a non-empty set with $0 \in U$ and $Y$ a normed linear space. We denote by $x(t)$ the unique solution of the initial-value problem:

$$\dot{x}(t) = f(t, T_r(t)x, u(t), d(t))$$
$$Y(t) = H(t, T_r(t)x) \quad (4.10)$$
$$x(t) \in \Re^n, Y(t) \in Y, d(t) \in D, u(t) \in U$$

with initial condition $T_r(t_0)x = x_0 \in C^0([-r,0];\Re^n)$, where $r > 0$ is a constant, $T_r(t)x := x(t+\theta)$; $\theta \in [-r,0]$ and the mappings $f : \Re_+ \times C^0([-r,0];\Re^n) \times U \times D \to \Re^n$, $H : \Re_+ \times C^0([-r,0];\Re^n) \to Y$ satisfy $f(t,0,0,d) = 0$, $H(t,0) = 0$ for all $(t,d) \in \Re_+ \times D$.

The following hypotheses will be imposed on systems of the form (4.10):

**(S1)** The mapping $(x,u,d) \to f(t,x,u,d)$ is continuous for each fixed $t \geq 0$ and there exists a symmetric positive definite matrix $P \in \Re^{n \times n}$ with the property that for every bounded $I \subseteq \Re_+$ and for every bounded $S \subset C^0([-r,0];\Re^n) \times U$, there exists a constant $L \geq 0$ such that:

$$(x(0) - y(0))' P(f(t,x,u,d) - f(t,y,u,d)) \leq L \max_{\tau \in [-r,0]} |x(\tau) - y(\tau)|^2 = L\|x - y\|_r^2$$

$\forall t \in I$, $\forall (x,u,y,u) \in S \times S$, $\forall d \in D$

**(S2)** There exist $a \in K_\infty$, $\gamma \in K^+$ such that $|f(t,x,u,d)| \leq \gamma(t)a(\|x\|_r + |u|)$ for all $(t,x,u,d) \in \Re_+ \times C^0([-r,0];\Re^n) \times U \times D$.

**(S3)** There exists a countable set $A \subset \Re_+$, which is either finite or $A = \{t_k \,;\, k = 1,...,\infty\}$ with $t_{k+1} > t_k > 0$ for all $k = 1,2,...$ and $\lim t_k = +\infty$, such that the mapping $(t,x,u,d) \in (\Re_+ \setminus A) \times C^0([-r,0];\Re^n) \times U \times D \to f(t,x,u,d)$ is continuous. Moreover, for each fixed $(t_0, x, u, d) \in \Re_+ \times C^0([-r,0];\Re^n) \times U \times D$, we have $\lim_{t \to t_0^+} f(t,x,u,d) = f(t_0,x,u,d)$.

**(S4)** The mapping $H : \Re_+ \times C^0([-r,0];\Re^n) \to Y$ is continuous.

For $x \in C^0([-r,0];\Re^n)$ we define $\|x\|_r := \max_{\theta \in [-r,0]} |x(\theta)|$. We will use the convention $C^0([0,0];\Re^n) = \Re^n$ and if $x \in C^0([0,0];\Re^n) = \Re^n$ we have $\|x\|_r = |x|$.

The class of functionals which are "almost Lipschitz on bounded sets" was introduced in [22,23] and is used extensively in the present work. For the sake of completeness we recall here the definition in [22,23].

**Definition 4.3:** *We say that a continuous functional $V : [-a,+\infty) \times C^0([-r,0];\Re^n) \to \Re_+$, $r > 0$, $a \geq 0$ is "almost Lipschitz on bounded sets", if there exist non-decreasing functions $M : \Re_+ \to \Re_+$, $P : \Re_+ \to \Re_+$, $G : \Re_+ \to [1,+\infty)$ such that for all $R \geq 0$, the following properties hold:*

**(P1)** *For every $x, y \in \{x \in C^0([-r,0];\Re^n)\,;\, \|x\|_r \leq R\}$, it holds that:*

$$|V(t,y) - V(t,x)| \leq M(R)\|y - x\|_r, \quad \forall t \in [-a, R]$$

**(P2)** *For every absolutely continuous function $x : [-r,0] \to \Re^n$ with $\|x\|_r \leq R$ and essentially bounded derivative, it holds that:*

$$|V(t+h,x) - V(t,x)| \leq hP(R)\left(1 + \operatorname*{ess\,sup}_{-r \leq \tau \leq 0} |\dot{x}(\tau)|\right), \text{ for all } t \in [-a,R] \text{ and } 0 \leq h \leq \frac{1}{G\left(R + \operatorname*{ess\,sup}_{-r \leq \tau \leq 0} |\dot{x}(\tau)|\right)}$$

*For the case $r = 0$, we say that a continuous functional $V : [-a,+\infty) \times C^0([-r,0];\Re^n) \to \Re_+$, is "almost Lipschitz on bounded sets", if $V : [-a,+\infty) \times \Re^n \to \Re_+$ is locally Lipschitz (notice that in this case $C^0([-r,0];\Re^n) = \Re^n$), i.e. for every compact $S \subset [-a,+\infty) \times \Re^n$ there exists $L \geq 0$ such that $|V(t,x) - V(\tau,y)| \leq L|t - \tau| + L|x - y|$ for all $(t,x) \in S$, $(\tau,y) \in S$.*

If the continuous functional $V:[-a,+\infty)\times C^0([-r,0];\Re^n) \to \Re_+$, is "almost Lipschitz on bounded sets" then we can define the derivative $V^0(t,x;v)$ in the following way (see also [22,23]) for $(t,x,v) \in \Re_+ \times C^0([-r,0];\Re^n) \times \Re^n$:

$$V^0(t,x;v) := \limsup_{h \to 0^+} \frac{V(t+h, E_h(x;v)) - V(t,x)}{h}$$

where $E_h(x;v)$ with $0 \le h < r$ denotes the following operator:

$$E_h(x;v) := \begin{cases} x(0) + (\theta + h)v & \text{for } -h < \theta \le 0 \\ x(\theta + h) & \text{for } -r \le \theta \le -h \end{cases} \quad (4.11a)$$

Particularly, for the case $r = 0$ we define

$$E_h(x;v) := x(0) + hv \quad (4.11b)$$

The following theorem provides sufficient Lyapunov-like conditions for the (U)IOS property. The gain functions of the IOS property can be determined **explicitly** in terms of the functions involved in the assumptions of the theorem. Its proof is provided in the Appendix.

**Theorem 4.4:** *Consider system (4.10) under hypotheses (S1-4) and suppose that there exist almost Lipschitz on bounded sets functionals $Q_i : [-r+r_i,+\infty) \times C^0([-r_i,0];\Re^n) \to \Re_+$ with $0 \le r_i \le r$ ($i=1,...,k$), $Q_0 : [-r+r_0,+\infty) \times C^0([-r_0,0];\Re^n) \to \Re_+$ with $0 \le r_0 \le r$, functions $a_1,a_2,a_3,a_4 \in K_\infty$ $\mu,\beta,\kappa \in K^+$, $\zeta \in N_1$, $g \in N_k$, $\gamma_{i,j} \in N_1$, $p_i \in N_1$, $i,j = 1,...,k$, positive definite functions $\rho_i \in C^0(\Re_+;\Re_+)$ ($i=1,...,k$) and a constant $\lambda \in (0,1)$ such that for all $(t,x,u) \in \Re_+ \times C^0([-r,0];\Re^n) \times U$ the following inequalities hold:*

$$a_1\left(\|H(t,x)\|_Y\right) \le \max_{i=1,...,k} V_i(t,x) \le a_2\left(\beta(t)\|x\|_r\right) \quad (4.12)$$

$$a_3\left(\mu(t)\|x\|_r\right) - g(V_1(t,x),...,V_k(t,x)) - \kappa(t) \le W(t,x) \le a_4\left(\beta(t)\|x\|_r\right) \quad (4.13)$$

$$\sup_{d \in D} Q_0^0(t, T_{r_0}(0)x; f(t,x,u,d)) \le -Q_0(t, T_{r_0}(0)x) + \lambda \max\left\{ \zeta(|u|), \max_{j=1,...,k} p_j(V_j(t,x)) \right\} \quad (4.14)$$

*where*

$$V_i(t,x) := \sup_{\theta \in [-r+r_i,0]} Q_i(t+\theta, T_{r_i}(\theta)x), \ i=1,...,k, \ W(t,x) := \sup_{\theta \in [-r+r_0,0]} Q_0(t+\theta, T_{r_0}(\theta)x) \quad (4.15)$$

*and for every $i=1,...,k$ and $(t,x,u) \in \Re_+ \times C^0([-r,0];\Re^n) \times U$ the following implication holds:*

*"If $\max\left\{ \zeta(|u|), \max_{j=1,...,k} \gamma_{i,j}(V_j(t,x)) \right\} \le Q_i(t, T_{r_i}(0)x)$ then $\sup_{d \in D} Q_i^0(t, T_{r_i}(0)x; f(t,x,u,d)) \le -\rho_i(Q_i(t, T_{r_i}(0)x))$"* (4.16)

*Further, suppose that the following set of small-gain conditions holds for $i=1,...,k$*

$$\gamma_{i,i}(s) < s, \ \forall s > 0 \quad (4.17)$$

*and for each $r = 2,...,k$:*

$$\left(\gamma_{i_1,i_2} \circ \gamma_{i_2,i_3} \circ ... \circ \gamma_{i_r,i_1}\right)(s) < s, \ \forall s > 0 \quad (4.18)$$

*for all $i_j \in \{1,...,k\}$, $i_j \ne i_l$ if $j \ne l$.*

*Then system (4.10) satisfies the IOS property with gain $\gamma = a_1^{-1} \circ \theta \in N_1$ from the input $u \in M_U$, where $\theta \in N_1$ is defined by (4.7), (4.8). Moreover, if $\beta \in K^+$ is bounded then system (4.10) satisfies the UIOS property with gain $\gamma = a_1^{-1} \circ \theta \in N_1$ from the input $u \in M_U$, where $\theta \in N_1$ is defined by (4.7), (4.8).*

When $H(t,x) = x$, setting $Q_0(t,x) \equiv 0$ in Theorem 4.4 we obtain the following corollary on the ISS of system (4.10).

**Corollary 4.5:** *Consider system (4.10) under hypotheses (S1-4) and suppose that there exists a family of almost Lipschitz on bounded sets functionals $Q_i : [-r+r_i, +\infty) \times C^0([-r_i, 0]; \Re^n) \to \Re_+$ with $0 \le r_i \le r$ ($i = 1,...,k$), functions $a_1, a_2 \in K_\infty$, $\beta \in K^+$, $\zeta \in N_1$, $\gamma_{i,j} \in N_1$, $i, j = 1,...,k$, and a family of positive definite functions $\rho_i \in C^0(\Re_+; \Re_+)$ ($i = 1,...,k$), such that for all $(t, x, u) \in \Re_+ \times C^0([-r, 0]; \Re^n) \times U$ the following inequality holds:*

$$a_1\left(\|x\|_r\right) \le \max_{i=1,...,k} V_i(t,x) \le a_2\left(\beta(t)\|x\|_r\right) \tag{4.19}$$

*where*

$$V_i(t,x) := \sup_{\theta \in [-r+r_i, 0]} Q_i(t+\theta, T_{r_i}(\theta)x), \; i = 1,...,k \tag{4.20}$$

*and implication (4.16) holds for every $i = 1,...,k$ and $(t, x, u) \in \Re_+ \times C^0([-r, 0]; \Re^n) \times U$. If, additionally, the set of small-gain conditions (4.17), (4.18) holds then system (4.10) satisfies the ISS property with gain $\gamma = a_1^{-1} \circ \theta \in N_1$ from the input $u \in M_U$, where $\theta \in N_1$ is defined by (4.7), (4.8). Moreover, if $\beta \in K^+$ is bounded then system (4.10) satisfies the UISS property with gain $\gamma = a_1^{-1} \circ \theta \in N_1$ from the input $u \in M_U$, where $\theta \in N_1$ is defined by (4.7), (4.8).*

**Remark 4.6:** It is of interest to note that some of the functionals $Q_i : [-r+r_i, +\infty) \times C^0([-r_i, 0]; \Re^n) \to \Re^+$ in Theorem 4.4 and Corollary 4.5 are allowed to be functions (case of $r_i = 0$). This reminds the case of Razumikhin functions, which are used frequently for the proof of stability properties of systems described by RFDEs (see [23,27,30,41]). Consequently, Theorem 4.4 and Corollary 4.5 allow the flexibility of using Lyapunov functionals with Razumikhin functions in order to prove desired stability properties.

**Remark 4.7:** It should be clear that the convention $C^0([0,0]; \Re^n) = \Re^n$ allows Theorem 4.4 and Corollary 4.5 to be used in the case of systems described by ODEs (case of $r = 0$). In this case Theorem 4.4 and Corollary 4.5 are generalizations of Theorem 4.1 and Corollary 4.2: the use of locally Lipschitz functions is allowed and discontinuities of the right-hand side of the differential equations with respect to time are allowed. In this case, the diagonal gain functions $\gamma_{i,i}$ for $i = 1,...,k$ play no role whatsoever and consequently can be set equal to zero. However, in the time-delay case (case of $r > 0$) the diagonal gain functions $\gamma_{i,i}$ for $i = 1,...,k$ play a significant role (see Example 5.1 below).

## 4.C. Sampled-data Systems

We consider switched systems, described in the following way: given a pair of sets $D \subseteq \Re^l$, $U \subseteq \Re^m$ with $0 \in U$, a positive function $h : \Re^n \times U \to (0, r]$, which is bounded by a certain constant $r > 0$ and a pair of vector fields $f : \Re^n \times \Re^n \times D \times U \times U \to \Re^n$, $H : \Re^n \to \Re^k$, we consider the switched system that produces for each $(t_0, x_0) \in \Re_+ \times \Re^n$ and for each triplet of measurable and locally bounded inputs $d : \Re_+ \to D$, $\tilde{d} : \Re_+ \to \Re_+$, $u : \Re_+ \to U$ the piecewise absolutely continuous function $t \to x(t) \in \Re^n$, via the following algorithm:

Step $i$:
1) Given $\tau_i$ and $x(\tau_i)$, calculate $\tau_{i+1}$ using the equation $\tau_{i+1} = \tau_i + \exp(-\tilde{d}(\tau_i))h(x(\tau_i), u(\tau_i))$,
2) Compute the state trajectory $x(t)$, $t \in [\tau_i, \tau_{i+1})$ as the solution of the differential equation $\dot{x}(t) = f(x(t), x(\tau_i), d(t), u(t), u(\tau_i))$,
3) Calculate $x(\tau_{i+1})$ using the equation $x(\tau_{i+1}) = \lim_{t \to \tau_{i+1}^-} x(t)$.

For $i = 0$ we take $\tau_0 = t_0$ and $x(\tau_0) = x_0$ (initial condition). Schematically, we write

$$\dot{x}(t) = f(x(t), x(\tau_i), d(t), u(t), u(\tau_i)) \quad , \quad t \in [\tau_i, \tau_{i+1})$$
$$\tau_0 = t_0, \tau_{i+1} = \tau_i + \exp(-\tilde{d}(\tau_i))h(x(\tau_i), u(\tau_i)), i = 0,1,... \qquad (4.21)$$
$$Y(t) = H(x(t))$$

with initial condition $x(t_0) = x_0$. Switched systems of the form (4.21) are called "sampled-data" systems (see also [29] and [39] for the case of state-dependent sampling period).

In the present work we study systems of the form (4.21) under the following hypotheses:

**(R1)** $f(x, x_0, d, u, u_0)$ is continuous with respect to $(x, d, u) \in \Re^n \times D \times U$ and such that for every bounded $S \subset \Re^n \times \Re^n \times U \times U$ there exists constant $L \geq 0$ such that

$$(x-y)'(f(x, x_0, d, u, u_0) - f(y, x_0, d, u, u_0)) \leq L|x-y|^2$$
$$\forall (x, x_0, u, u_0, d) \in S \times D, \forall (y, x_0, u, u_0, d) \in S \times D$$

**(R2)** There exists a function $a \in K_\infty$ such that

$$|f(x, x_0, d, u, u_0)| \leq a(|x| + |x_0| + |u| + |u_0|), \forall (u, u_0, d, x, x_0) \in U \times U \times D \times \Re^n \times \Re^n$$

**(R3)** $H: \Re^n \to \Re^k$ is a continuous map with $H(0) = 0$.

**(R4)** The function $h: \Re^n \times U \to (0, r]$ is a positive, continuous and bounded function.

The following theorem provides sufficient Lyapunov-like conditions for the (U)IOS property. The gain functions of the IOS property can be determined **explicitly** in terms of the functions involved in the assumptions of the theorem. Its proof is provided in the Appendix.

**Theorem 4.8 (Vector Lyapunov Function Characterization of UIOS):**

*Consider system (4.21) under hypotheses (R1-4) and suppose that there exist nonnegative functions $V_i \in C^1(\Re^n; \Re_+)$ ($i = 1,...,k$), $Q \in C^1(\Re^n; \Re_+)$, $a_1, a_2, a_3, a_4 \in K_\infty$, $\zeta \in N_1$, $g \in N_k, \gamma_{i,j} \in N_1$, $p_i \in N_1$, $i, j = 1,...,k$, constants $\mu, \kappa \geq 0$, $\lambda \in (0,1)$ and positive definite functions $\rho_i \in C^0(\Re_+; \Re_+)$ ($i = 1,...,k$), such that the following inequalities hold for all $(x, x_0, u, u_0) \in \Re^n \times \Re^n \times U \times U$:*

$$a_1(|H(x)|) \leq \max_{i=1,...,k} V_i(x) \leq a_2(|x|) \qquad (4.22)$$

$$a_3(|x|) - g(V_1(x),...,V_k(x)) - \kappa \leq Q(x) \leq a_4(|x|) \qquad (4.23)$$

$$\sup_{d \in D} \nabla Q(x) f(x, x_0, d, u, u_0) \leq \mu Q(x) + \lambda \max\left\{\zeta(|u|), \zeta(|u_0|), \max_{j=1,...,k} p_j(V_j(x)), \max_{j=1,...,k} p_j(V_j(x_0))\right\} \qquad (4.24)$$

*and for every $i = 1,...,k$ and $(x, u, u_0) \in \Re^n \times U \times U$ the following implication holds:*

"If $\max\left\{\zeta(|u|), \zeta(|u_0|), \max_{j=1,...,k} \gamma_{i,j}(V_j(x))\right\} \leq V_i(x)$ and $x_0 \in A_i(h(x_0, u_0), x)$ then
$$\sup_{d \in D} \nabla V_i(x) f(x, x_0, d, u, u_0) \leq -\rho_i(V_i(x))" \qquad (4.25)$$

*where the family of set-valued maps $\Re_+ \times \Re^n \ni (T, x) \to A_i(T, x) \subseteq \Re^n$ ($i = 1,...,k$) is defined by*

$$A_i(T,x) = \bigcup_{0 \le s \le T} \left\{ \begin{array}{l} x_0 \in \Re^n : \exists (d,u) \in M_D \times M_U \text{ with } \phi(s, x_0; d, u) = x, \\ \zeta(|u(t)|) \le V_i(x), \gamma_{i,j}(V_j(\phi(t, x_0; d, u))) \le V_i(x) \\ \text{for all } t \in [0, s] \text{ and } j = 1, ..., k \end{array} \right\} \quad (4.26)$$

and $\phi(t, x_0; d, u)$ denotes the solution of $\dot{x}(t) = f(x(t), x_0, d(t), u(t), u(0))$ with initial condition $x(0) = x_0$ corresponding to $(d, u) \in M_D \times M_U$.

*Furthermore, if the set of small-gain conditions (4.17), (4.18) holds then system (4.21) satisfies the UIOS property with gain $\gamma = a_1^{-1} \circ \theta \in N_1$ from the input $u \in M_U$ and zero gain from the input $\tilde{d} \in M_{\Re_+}$, where $\theta \in N_1$ is defined by (4.7), (4.8).*

For the ISS case where $H(t, x) = x$, one can set $Q(x) \equiv 0$ in Theorem 4.8 and obtain the following corollary.

**Corollary 4.9 (Vector Lyapunov Function Characterization of UISS):**

*Consider system (4.21) under hypotheses (R1-4) and suppose that there exists a family of functions $V_i \in C^1(\Re^n; \Re_+)$ ($i = 1, ..., k$), functions $a_1, a_2 \in K_\infty$, $\zeta \in N_1$, $\gamma_{i,j} \in N_1$, $i, j = 1, ..., k$, and a family of positive definite functions $\rho_i \in C^0(\Re_+; \Re_+)$ ($i = 1, ..., k$), such that the following inequality holds for all $x \in \Re^n$:*

$$a_1(|x|) \le \max_{i=1,...,k} V_i(x) \le a_2(|x|) \quad (4.27)$$

*and implication (4.25) holds for every $i = 1, ..., k$ and $(x, u, u_0) \in \Re^n \times U \times U$, where the family of set-valued maps $\Re_+ \times \Re^n \ni (T, x) \to A_i(T, x) \subseteq \Re^n$ ($i = 1, ..., k$) is defined by (4.26) and $\phi(t, x_0; d, u)$ denotes the solution of $\dot{x}(t) = f(x(t), x_0, d(t), u(t), u(0))$ with initial condition $x(0) = x_0$ corresponding to $(d, u) \in M_D \times M_U$.*

*Under the set of small-gain conditions (4.17), (4.18), system (4.21) satisfies the UISS property with gain $\gamma = a_1^{-1} \circ \theta \in N_1$ from the input $u \in M_U$ and zero gain from the input $\tilde{d} \in M_{\Re_+}$, where $\theta \in N_1$ is defined by (4.7), (4.8).*

**Remark 4.10:** It is worth noting that Theorem 3.1 and Corollary 3.3 in [24] are special cases of Theorem 4.8 and Corollary 4.9 with $\gamma_{i,j}(s) = a(s)$ for all $i, j = 1, ..., k$ and $Q(x) \equiv 0$, where $a \in N_1$ with $a(s) < s$ for $s > 0$. Moreover, it is assumed in [24] that there exist a constant $R \ge 0$ and a function $p \in K_\infty$ such that $|x| \le R + p(|H(x)|)$. In the present work, such a hypothesis is not needed.

The interpretation of the family of set-valued maps $\Re_+ \times \Re^n \ni (T, x) \to A_i(T, x) \subseteq \Re^n$ ($i = 1, ..., k$), defined in (4.26) is the following (the same with [24]): each $A_i(T, x) \subseteq \Re^n$ is the set of all states $x_0 \in \Re^n$ so that the solution of $\dot{x}(t) = f(x(t), x_0, d(t), u(t), u(0))$ with initial condition $x(0) = x_0$ can be controlled to $x \in \Re^n$ in time $s$ less or equal than $T$ by means of appropriate inputs $(d, u) \in M_D \times M_U$ that satisfy $\zeta\left(\sup_{t \in [0,s]} |u(t)|\right) \le V_i(x)$ and such that the trajectory of the solution satisfies the constraint $\max_{j=1,...,k} \sup_{t \in [0,s]} \gamma_{i,j}(V_j(x(t))) \le V_i(x)$. In general it is very difficult to obtain an accurate description of the set-valued maps $\Re_+ \times \Re^n \ni (T, x) \to A_i(T, x) \subseteq \Re^n$ defined by (4.26). However, for every $g \in C^1(\Re^n; \Re)$, we have:

$$A_i(T, x) \subseteq B_i^g(T, x) = \left\{ x_0 \in \Re^n : |g(x_0) - g(x)| \le T b_i^g(x) \right\}, \quad \forall (T, x) \in \Re_+ \times \Re^n$$

where

$$b_i^g(x) := \max\left\{ |\nabla g(\xi) f(\xi, x_0, d, u, u_0)| : d \in D, \zeta(\max\{|u|, |u_0|\}) \le V_i(x), \max_{j=1,...,k} \gamma_{i,j}(\max\{V_j(\xi), V_j(x_0)\}) \le V_i(x) \right\} < +\infty$$

and $V_i \in C^1(\Re^n; \Re_+)$ ($i = 1, ..., k$) are the functions involved in hypotheses of Theorem 4.8.

## 5. Examples and Applications

**Example 5.1:** Consider the time-delay system:
$$\dot{x}_i(t) = -a_i x_i(t) + g_i(d(t), T_r(t)x) \quad , \quad i = 1,...,n \tag{5.1}$$

where $d(t) \in D \subseteq \Re^m$, $a_i > 0$ ($i = 1,...,n$) and $g_i : D \times C^0([-r,0]; \Re^n) \to \Re$ ($i = 1,...,n$) are continuous mappings with

$$\sup_{d \in D} |g_i(d,x)| \leq \max_{j=1,...,n} c_{i,j} \|x_j\|_r \tag{5.2}$$

for certain constants $c_{i,j} \geq 0$ ($i, j = 1,...,n$). We next show that $0 \in C^0([-r,0]; \Re^n)$ is RGAS for (5.1) if $c_{i,i} < a_i$ for all $i = 1,...,n$ and the following small-gain conditions hold for each $r = 2,...,n$:

$$c_{i_1,i_2} c_{i_2,i_3} ... c_{i_r,i_1} < a_{i_1} a_{i_2} ... a_{i_r} \tag{5.3}$$

for all $i_j \in \{1,...,n\}$, $i_j \neq i_k$ if $j \neq k$.

First, we notice that hypotheses (S1-4) hold for system (5.1) under hypothesis (5.2) with output $H(t,x) := x \in C^0([-r,0]; \Re^n)$. Define the family of functions $Q_i(x) = \frac{1}{2} x_i^2(0)$ and $V_i(x) := \sup_{\theta \in [-r,0]} Q_i(x(\theta)) = \frac{1}{2} \|x_i\|_r^2$ ($i = 1,...,n$) for $x \in C^0([-r,0]; \Re^n)$. These mappings satisfy inequality (4.19) and definition (4.20) with $a_1(s) := \frac{1}{2n} s^2$, $a_2(s) := \frac{1}{2} s^2$, $\beta(t) \equiv 1$ and $r_i := 0$, $i = 1,...,n$. Let $\lambda \in (0,1)$ and notice that implication (4.16) holds with $\gamma_{i,j}(s) := \frac{c_{i,j}^2}{\lambda^2 a_i^2} s$ and $\rho_i(s) := 2(1-\lambda) a_i s$. Condition (5.3) and the fact that $c_{i,i} < a_i$ for all $i = 1,...,n$ implies that the small-gain conditions (4.17), (4.18) hold for $\lambda \in (0,1)$ sufficiently close to 1. We conclude from Corollary 4.5 that $0 \in C^0([-r,0]; \Re^n)$ is RGAS for (5.1).

It is important to notice that the conditions on the diagonal terms cannot be avoided in general. Such a situation occurs for example when $\dot{x}_1(t) = -a_1 x_1(t) + c_{1,1} x_1(t-r) + c_{1,2} x_2(t-r)$ with $c_{1,1} \neq 0$. ◁

**Example 5.2:** Consider the following biochemical control circuit model:
$$\begin{aligned} \dot{X}_1(t) &= g(X_n(t-\tau_n)) - a_1 X_1(t) \\ \dot{X}_i(t) &= X_{i-1}(t-\tau_{i-1}) - a_i X_i(t) \quad , \quad i = 2,...,n \\ X(t) &= (X_1(t),...,X_n(t))' \in \Re_+^n \end{aligned} \tag{5.4}$$

where $a_i > 0$, $\tau_i \geq 0$ ($i = 1,...,n$) are constants and $g \in C^1(\Re_+; \Re_+)$ is a function with $g(X) > 0$ for all $X > 0$. This model has been studied in [32] (see pages 58-60 and 93-94). In this book it is further assumed that $g \in C^1(\Re_+; \Re_+)$ is bounded and strictly increasing (a typical choice for $g \in C^1(\Re_+; \Re_+)$ is $g(X) = \frac{X^p}{1 + X^p}$ with $p$ being a positive integer). It is shown that if there is one equilibrium point for (5.4) then it attracts all solutions. If there are two equilibrium points then all solutions are attracted to these points. Here we study (5.4) under the following assumption:

**(H)** There exist $X_n^* > 0$, $K > 0$ and $\lambda \in (0,1)$ with $aX_n^* = g(X_n^*)$ and such that
$$\frac{K + X_n^*}{K + X} X \leq a^{-1} g(X) \leq X_n^* + \lambda |X - X_n^*|, \text{ for all } X \geq 0 \tag{5.5}$$

where $a = \prod_{j=1}^n a_j$.

Using Small-Gain Analysis we are in a position to prove:

> "Consider system (5.4) under hypothesis (H) and let $r := \max_{i=1,...,n} \tau_i$. Then for every $X_0 \in C^0([-r,0]; \text{int}(\Re_+^n))$ the solution of (5.4) with initial condition $T_r(0)X = X_0$ satisfies $\lim_{t \to +\infty} X(t) = X^*$, where $X^* = (X_1^*,...,X_n^*)' \in \text{int}(\Re_+^n)$ with $\left(\prod_{j=1}^{i} a_j\right) X_i^* = g(X_n^*)$, for $i = 1,...,n-1$."

It should be clear that in contrast to the analysis performed in [32] for (5.4) (based on the monotone dynamical system theory) we do not assume that $g \in C^1(\Re_+; \Re_+)$ is bounded or strictly increasing, Moreover, even if there are two equilibrium points (notice that (5.5) allows $g(0) = 0$ and therefore $0 \in \Re_+^n$ can be an equilibrium point), we prove almost global convergence to the non-trivial equilibrium.

A typical analysis of the equilibrium points of (5.4) under hypothesis (H) shows that there exists an equilibrium point $X^* \in \text{int}(\Re_+^n)$ satisfying:

$$\left(\prod_{j=1}^{i} a_j\right) X_i^* = g(X_n^*), \quad i = 1,...,n \tag{5.6}$$

In order to be able to study solutions of (5.4) evolving in $\text{int}(\Re_+^n)$ we consider the following transformation:

$$X_i = X_i^* \exp(x_i), \quad i = 1,...,n \tag{5.7}$$

Therefore system (5.4) under transformation (5.7) is expressed by the following set of differential equations:

$$\dot{x}_1 = a_1 \left( \frac{g(X_n^* \exp(x_n(t - \tau_n)))}{g(X_n^*)} \exp(-x_1(t)) - 1 \right) \tag{5.8a}$$

$$\dot{x}_i(t) = a_i \left( \exp(x_{i-1}(t - \tau_{i-1}) - x_i(t)) - 1 \right), \quad i = 2,...,n$$
$$x(t) = (x_1(t),...,x_n(t))' \in \Re^n \tag{5.8b}$$

First, we notice that hypotheses (S1-4) hold for system (5.8) under hypothesis (H) with output $H(t,x) := x \in C^0([-r,0]; \Re^n)$ and that $0 \in C^0([-r,0]; \Re^n)$ is an equilibrium point for (5.8). Define the family of functions $Q_i(x) = \frac{1}{2} x_i^2(0)$ and $V_i(x) := \sup_{\theta \in [-r,0]} Q_i(x(\theta)) = \frac{1}{2} \|x_i\|_r^2$ ($i = 1,...,n$) for $x \in C^0([-r,0]; \Re^n)$. These mappings satisfy inequality (4.19) and definition (4.20) with $a_1(s) := \frac{1}{2n} s^2$, $a_2(s) := \frac{1}{2} s^2$, $\beta(t) \equiv 1$ and $r_i := 0$, $i = 1,...,n$.

We define $\gamma_{1,j}(s) \equiv 0$ for $j \neq n$ and $\gamma_{1,n}(s) := \frac{1}{2} \left[\log\left(1 + \theta\left(\exp(\sqrt{2s}) - 1\right)\right)\right]^2$, where $\theta \in \left(\max\left\{\frac{b}{b+1}, \lambda\right\}, 1\right)$, $\lambda \in (0,1)$ being the constant involved in hypothesis (H) and $b := \frac{K}{X_n^*}$. Notice that

$$Q_1^0\left(x_1(0); a_1\left(\frac{g(X_n^* \exp(x_n(-\tau_n)))}{g(X_n^*)} \exp(-x_1(0)) - 1\right)\right) = a_1 x_1(0) \left(\frac{g(X_n^* \exp(x_n(-\tau_n)))}{g(X_n^*)} \exp(-x_1(0)) - 1\right)$$

We consider the following cases:

1) $x_1(0) < 0$. In this case the left hand side inequality (5.5) implies that
$\frac{g(X_n^* \exp(x_n(-\tau_n)))}{g(X_n^*)} \geq \frac{b+1}{b+\exp(x_n(-\tau_n))} \exp(x_n(-\tau_n)) \geq \frac{b+1}{b+\exp(-|x_n(-\tau_n)|)} \exp(-|x_n(-\tau_n)|)$, where $b := \frac{K}{X_n^*}$. The inequality $\gamma_{1,n}(V_n(x)) \leq Q_1(x_1(0))$ implies $\ln(1+\theta(\exp(|x_n(-\tau_n)|)-1)) \leq -x_1(0)$, which combined with the previous inequalities gives:

$$Q_1^0\left(x_1(0); a_1\left(\frac{g(X_n^* \exp(x_n(-\tau_n)))}{g(X_n^*)} \exp(-x_1(0)) - 1\right)\right) \leq a_1 x_1(0) \frac{(b+1-b\theta^{-1})(\exp(-x_1(0))-1)}{b+1+b\theta^{-1}(\exp(-x_1(0))-1)} \quad (5.9)$$

2) $x_1(0) \geq 0$. In this case the right hand side inequality (5.5) implies that
$\frac{g(X_n^* \exp(x_n(-\tau_n)))}{g(X_n^*)} \leq 1 + \lambda|\exp(x_n(-\tau_n))-1| \leq 1 + \lambda(\exp(|x_n(-\tau_n)|)-1)$. The inequality $\gamma_{1,n}(V_n(x)) \leq Q_1(x_1(0))$ implies $\ln(1+\theta(\exp(|x_n(-\tau_n)|)-1)) \leq x_1(0)$, which combined with the previous inequalities gives:

$$Q_1^0\left(x_1(0); a_1\left(\frac{g(X_n^* \exp(x_n(-\tau_n)))}{g(X_n^*)} \exp(-x_1(0)) - 1\right)\right) \leq a_1 x_1(0)(\lambda\theta^{-1}-1)(1-\exp(-x_1(0))) \quad (5.10)$$

Combining the two cases we obtain from (5.9) and (5.10) that the following implication holds:

$$\gamma_{1,n}(V_n(x)) \leq Q_1(x_1(0)) \Rightarrow Q_1^0\left(x_1(0); a_1\left(\frac{g(X_n^* \exp(x_n(-\tau_n)))}{g(X_n^*)} \exp(-x_1(0)) - 1\right)\right) \leq -\rho_1(Q_1(x_1(0))) \quad (5.11)$$

with $\rho_1(s) := a_1\sqrt{2s} \min\left\{(1-\lambda\theta^{-1})(1-\exp(-\sqrt{2s})), \frac{(b+1-b\theta^{-1})(\exp(\sqrt{2s})-1)}{b+1+b\theta^{-1}(\exp(\sqrt{2s})-1)}\right\}$.

We next define for $i=2,\ldots,n$, $\gamma_{i,j}(s) \equiv 0$ for $j \neq i-1$ and $\gamma_{i,i-1}(s) := \frac{1}{2}[\log(1+\mu(\exp(\sqrt{2s})-1))]^2$, where $\mu > 1$ is to be selected. Working in a similar way as above we obtain for $i=2,\ldots,n$:

$$\gamma_{i,i-1}(V_i(x)) \leq Q_i(x_i(0)) \Rightarrow Q_i^0(x_i(0); a_i(\exp(x_{i-1}(-\tau_{i-1})-x_i(0))-1)) \leq -\rho_i(Q_i(x_i(0))) \quad (5.12)$$

with $\rho_i(s) := (1-\mu^{-1})a_i\sqrt{2s}\frac{(1-\exp(-\sqrt{2s}))}{1+\mu^{-1}(\exp(\sqrt{2s})-1)}$ for $i=2,\ldots,n$.

Therefore, we conclude from (5.11) and (5.12) that implication (4.16) holds.

Finally, we check the small-gain conditions. Exploiting the previous definitions of the functions $\gamma_{i,j}(s)$, $i,j=1,\ldots,n$, we conclude that the small-gain conditions (4.17), (4.18) hold if and only if $(\gamma_{n,n-1} \circ \gamma_{n-1,n-2} \circ \ldots \circ \gamma_{1,n})(s) < s$ for all $s > 0$. Since

$$(\gamma_{n,n-1} \circ \gamma_{n-1,n-2} \circ \ldots \circ \gamma_{1,n})(s) = \frac{1}{2}[\log(1+\mu^{n-1}\theta(\exp(\sqrt{2s})-1))]^2$$

the small-gain conditions (4.17), (4.18) hold with $\mu \in \left(1, \theta^{-\frac{1}{n-1}}\right)$.

Thus, Corollary 4.5 implies that $0 \in C^0([-r,0]; \Re^n)$ is Uniformly Globally Asymptotically Stable for system (5.8). Taking into account transformation (5.7), this implies that for every $X_0 \in C^0([-r,0]; \text{int}(\Re_+^n))$ the solution of (5.4) with initial condition $T_r(0)X = X_0$ satisfies $\lim_{t \to +\infty} X(t) = X^*$, where $X^* = (X_1^*,...,X_n^*)' \in \text{int}(\Re_+^n)$ with $\left(\prod_{j=1}^{i} a_j\right) X_i^* = g(X_n^*)$, for $i = 1,...,n-1$. ◁

## 6. Conclusions

A novel Small-Gain Theorem is presented, which leads to vector Lyapunov characterizations of the (uniform and non-uniform) IOS property for various important classes of nonlinear control systems. The results presented in this work generalize many recent small-gain results in the literature and allow the explicit computation of the gain function of the overall system. Moreover, since the gain map $\Gamma : \Re_+^n \to \Re_+^n$ is allowed to contain diagonal terms, the obtained results have direct applications to time-delay systems. Examples have demonstrated the effectiveness of the vector small-gain methodology to large-scale time-delay systems, such as those encountered in biotechnology.

Our future work will be directed at applications of the vector small-gain theorem to the nonlinear feedback design issue for various classes of nonlinear control systems. Another interesting topic for future research is to study the internal and external stability properties for coupled systems involving integral input-to-state stable (iISS, a weaker notion than ISS [36]) subsystems from a viewpoint of vector small-gain. Some preliminary results are reported upon in [8] for interconnected systems consisting of two ISS and/or iISS subsystems.

## Appendix-Proofs

**Proof of Proposition 2.7:** We prove implications (iv) $\Rightarrow$ (i), (iii) $\Rightarrow$ (ii) and (ii) $\Rightarrow$ (iv), since the implication (i) $\Rightarrow$ (iii) is a consequence of Proposition 2.1.

(iv) $\Rightarrow$ (i): Clearly, since $\Gamma^{(k)}(x) \leq Q(x)$ for all $k \geq 1$, $x \in \Re_+^n$, we have $|\Gamma^{(k)}(x)| \leq |Q(x)|$ for all $k \geq 1$, $x \in \Re_+^n$ (notice that $\Gamma^{(k)}(x) \in \Re_+^n$). Continuity of the mapping $Q(x) = MAX\{x, \Gamma(x), \Gamma^{(2)}(x),...,\Gamma^{(n-1)}(x)\}$ (which is a direct consequence of continuity of the mapping $\Gamma(x)$) implies that for every $\varepsilon > 0$ there exists $\delta > 0$ such that $|x| \leq \delta$, $x \in \Re_+^n$ implies $|Q(x)| \leq \varepsilon$ (notice that $Q(0) = 0$). This implies stability.

Since $\Gamma : \Re_+^n \to \Re_+^n$ is MAX-preserving we have $\Gamma(Q(x)) = MAX\{\Gamma(x), \Gamma^{(2)}(x),...,\Gamma^{(n)}(x)\}$ for all $x \in \Re_+^n$. Moreover, since $\Gamma^{(k)}(x) \leq Q(x)$ for all $k \geq 1$, $x \in \Re_+^n$ it follows that $\Gamma(Q(x)) \leq Q(x)$ for all $x \in \Re_+^n$. Lemma 2.2 in conjunction with the fact that (iii) holds and $x \leq Q(x)$ for all $x \in \Re_+^n$ implies that $\lim_{k \to \infty} \Gamma^{(k)}(x) = 0$ for all $x \in \Re_+^n$.

(iii) $\Rightarrow$ (ii): If there exist $s > 0$ and some integer $i = 1,...,n$ such that $\gamma_{i,i}(s) \geq s$ then the non-zero vector $x \in \Re_+^n$ with $x_i = s$ and $x_j = 0$ for $j \neq i$ will violate (iii). Consequently, $\gamma_{i,i}(s) < s$, for all $s > 0$, $i = 1,...,n$.

Next suppose that $n > 1$. Suppose that there exist some $s > 0$, $r \in \{2,...,n\}$, indices $i_j \in \{1,...,n\}$, $j = 1,...,r$ with $i_j \neq i_k$ if $j \neq k$ such that $(\gamma_{i_1,i_2} \circ \gamma_{i_2,i_3} \circ ... \circ \gamma_{i_r,i_1})(s) \geq s$. Without loss of generality we may assume that $i_j = j$, for $j = 1,...,r$ and consequently $(\gamma_{1,2} \circ \gamma_{2,3} \circ ... \circ \gamma_{r,1})(s) \geq s$. The non-zero vector $x \in \Re_+^n$ with $x_1 = s$, $x_j = (\gamma_{j,j+1} \circ \gamma_{j+1,j+2} \circ ... \circ \gamma_{r,1})(s)$ for $j = 2,...,r$ and $x_j = 0$ for $j > r$ satisfies $\Gamma(x) \geq x$ and consequently hypothesis (iii) is violated. Therefore (ii) must hold.

(ii) $\Rightarrow$ (iv) The proof of this implication is a direct consequence of the fact that

$$\Gamma_i^{(k)}(x) = \max\left\{(\gamma_{i,j_1} \circ \gamma_{j_1,j_2} ... \circ \gamma_{j_{k-1},j_k})(x_{j_k}) : (j_1,...,j_k) \in \{1,...,n\}^k\right\}$$

for all $k \geq 1$, $x \in \Re_+^n$ and $i = 1,...,n$. Using (ii) it may be shown that $\Gamma^{(n)}(x) \leq Q(x) = MAX\{x, \Gamma(x), \Gamma^{(2)}(x),..., \Gamma^{(n-1)}(x)\}$ for all $x \in \Re_+^n$. Since $\Gamma : \Re_+^n \to \Re_+^n$ is MAX-preserving we have $\Gamma(Q(x)) = MAX\{\Gamma(x), \Gamma^{(2)}(x),..., \Gamma^{(n)}(x)\}$ for all $x \in \Re_+^n$. As a result, we obtain $\Gamma(Q(x)) \leq Q(x)$ for all $x \in \Re_+^n$. By induction, it follows that $\Gamma^{(k)}(Q(x)) \leq Q(x)$ for all $k \geq 1$, $x \in \Re_+^n$. Since $x \leq Q(x)$, we obtain $\Gamma^{(k)}(x) \leq Q(x)$ for all $k \geq 1$, $x \in \Re_+^n$.

The fact that implication (iii) holds is shown by contradiction. Suppose that there exists a non-zero $x \in \Re_+^n$ with $\Gamma(x) \geq x$. Consequently, for every $i \in \{1,...,n\}$ there exists $p(i) \in \{1,...,n\}$ with $\gamma_{i,p(i)}(x_{p(i)}) \geq x_i$. With these inequalities in mind, there is at least one $i \in \{1,...,n\}$ with $x_i > 0$, and a closed cycle $(i, j_1,..., j_r, i)$ such that $\left(\gamma_{i,j_1} \circ \gamma_{j_1,j_2} \circ ... \circ \gamma_{j_r,i}\right)(x_i) \geq x_i$, which contradicts (ii). Therefore the implication (ii) => (iv) holds.

The proof is thus completed.  ◁

**Proof of Theorem 3.1:** The proof consists of two steps:

Step 1: We show that $\Sigma$ is RFC from the input $u \in M_U$ and that for every $(t_0, x_0, u, d) \in \Re_+ \times \mathcal{X} \times M_U \times M_D$ the following inequality holds for all $t \geq t_0$:

$$V(t) \leq MAX\left\{ Q(\mathbf{1}\sigma(L(t_0),0)), Q\left(\mathbf{1}\zeta\left(\|u(\tau)\|_{\mathcal{U}}\big]_{[t_0,t]}\right)\right)\right\} \tag{A.1}$$

Therefore by virtue of (A.1), (3.3), properties P1 and P2 of Lemma 2.7 in [20] hold for system $\Sigma$ with $V = V_i$ and $\gamma = G_i$ ($i = 1,...,n$). Moreover, if $\beta \in K^+$ is bounded then (3.3) implies that properties P1 and P2 of Lemma 2.8 in [20] hold for system $\Sigma$ with $V = V_i$ and $\gamma = G_i$ ($i = 1,...,n$).

Step 2: We prove the following claim.

**Claim:** *For every $\varepsilon > 0$, $k \in Z_+$, $R, T \geq 0$ there exists $\tau_k(\varepsilon, R, T) \geq 0$ such that for every $(t_0, x_0, u, d) \in \Re_+ \times \mathcal{X} \times M_U \times M_D$ with $t_0 \in [0,T]$ and $\|x_0\|_{\mathcal{X}} \leq R$ the following inequality holds:*

$$V(t) \leq MAX\left\{ Q(P\varepsilon), \Gamma^{(k)}(Q(\mathbf{1}\sigma(L(t_0),0))), G\left(\|u(\tau)\|_{\mathcal{U}}\big]_{[t_0,t]}\right)\right\}, \text{ for all } t \geq t_0 + \tau_k \tag{A.2}$$

*Moreover, if $\beta, c \in K^+$ are bounded then for every $\varepsilon > 0$, $k \in Z_+$, $R \geq 0$ there exists $\tau_k(\varepsilon, R) \geq 0$ such that for every $(t_0, x_0, u, d) \in \Re_+ \times \mathcal{X} \times M_U \times M_D$ with $\|x_0\|_{\mathcal{X}} \leq R$ inequality (A.2) holds.*

Notice that hypothesis (H2) and inequality (3.3) guarantees the existence of $k(\varepsilon, T, R) \in Z_+$ such that $Q(P\varepsilon) \geq \Gamma^{(l)}\left(Q\left(\mathbf{1}\sigma\left(b\left(R \max_{0 \leq t \leq T} \beta(t)\right),0\right)\right)\right)$ for all $l \geq k$. If $\beta \in K^+$ is bounded then $k$ is independent of $T$. Therefore by virtue of (A.2), property P3 of Lemma 2.7 in [20] holds for system $\Sigma$ with $V = V_i$ and $\gamma = G_i$ ($i = 1,...,n$). Moreover, if $\beta, c \in K^+$ are bounded then (A.2) implies that property P3 of Lemma 2.8 in [20] hold for system $\Sigma$ with $V = V_i$ and $\gamma = G_i$ ($i = 1,...,n$).

The proof of Theorem 3.1 is thus completed with the help of Lemma 2.7 (or Lemma 2.8) in [20].

Now, we return to establish RFC as claimed in Step 1.

Step 1:

Let $(t_0, x_0, u, d) \in \Re_+ \times X \times M_U \times M_D$. Inequality (3.1) implies for all $t \in [t_0, t_{\max})$

$$[V]_{[t_0,t]} \leq MAX\left\{ \mathbf{1}\sigma(L(t_0),0), \Gamma([V]_{[t_0,t]}), \mathbf{1}\zeta\left( \left[\|u(\tau)\|_U\right]_{[t_0,t]}\right) \right\} \tag{A.3}$$

Proposition 2.9 in conjunction with (A.3) implies (A.1) for all $t \in [t_0, t_{\max})$.

We show next that $\Sigma$ is RFC from the input $u \in M_U$ by contradiction. Suppose that $t_{\max} < +\infty$. Then by virtue of the BIC property for every $M > 0$ there exists $t \in [t_0, t_{\max})$ with $\|\phi(t, t_0, x_0, u, d)\|_X > M$. On the other hand estimate (A.1) in conjunction with the hypothesis $t_{\max} < +\infty$ shows that there exists $M_1 \geq 0$ such that $\sup_{t_0 \leq \tau < t_{\max}} |V(\tau)| \leq M_1$. The fact that $V(t)$ is bounded in conjunction with estimate (3.2) implies that there exist constants $M_2 \geq 0$ such that $\sup_{t_0 \leq \tau < t_{\max}} L(\tau) \leq M_2$. It follows from (3.3) and inequality $\mu(t)\|\phi(t, t_0, x_0, u, d)\|_X \leq b(L(t) + g(V(t)) + \kappa(t))$ that the transition map of $\Sigma$, i.e., $\phi(t, t_0, x_0, u, d)$, is bounded on $[t_0, t_{\max})$ and this contradicts the requirement that for every $M > 0$ there exists $t \in [t_0, t_{\max})$ with $\|\phi(t, t_0, x_0, u, d)\|_X > M$. Hence, we must have $t_{\max} = +\infty$.

Therefore, we conclude that $\Sigma$ is RFC from the input $u \in M_U$ and that (A.1) holds for all $t \geq t_0$.

Step 2: Proof of the Claim

The proof of the claim will be made by induction on $k \in Z_+$.

First we show inequality (A.2) for $k = 1$.

Let arbitrary $\varepsilon > 0$, $R, T \geq 0$, $(t_0, x_0, u, d) \in \Re_+ \times X \times M_U \times M_D$ with $t_0 \in [0, T]$ and $\|x_0\|_X \leq R$. Inequality (3.1) in conjunction with inequality (A.1) give for $t \geq t_0$:

$$V(t) \leq MAX\left\{ \mathbf{1}\sigma(L(t_0), t-t_0), \Gamma(Q(\mathbf{1}\sigma(L(t_0),0))), \Gamma\left(Q\left(\mathbf{1}\gamma\left(\left[\|u(\tau)\|_U\right]_{[t_0,t]}\right)\right)\right), \mathbf{1}\zeta\left(\left[\|u(\tau)\|_U\right]_{[t_0,t]}\right) \right\} \tag{A.4}$$

Since $\Gamma(Q(x)) \leq Q(x)$ and $Q(x) \geq x$ for all $x \in \Re_+^n$, inequality (A.4) implies for all $t \geq t_0$:

$$V(t) \leq MAX\left\{ \mathbf{1}\sigma(L(t_0), t-t_0), \Gamma(Q(\mathbf{1}\sigma(L(t_0),0))), Q\left(\mathbf{1}\zeta\left(\left[\|u(\tau)\|_U\right]_{[t_0,t]}\right)\right) \right\} \tag{A.5}$$

Similarly inequality (3.2) in conjunction with inequality (A.1) give for $t \geq t_0$:

$$L(t) \leq \max\left\{ \nu(t-t_0), c(t_0), a(\|x_0\|_X), p(Q(\mathbf{1}\sigma(L(t_0),0))), p\left(Q\left(\mathbf{1}\zeta\left(\left[\|u(\tau)\|_U\right]_{[t_0,t]}\right)\right)\right), p^u\left(\left[\|u(\tau)\|_U\right]_{[t_0,t]}\right) \right\} \tag{A.6}$$

Notice that (3.3) implies $L(t_0) \leq b(\beta(t_0)\|x_0\|_X) \leq b\left(R \max_{0 \leq t \leq T} \beta(t)\right)$. Using the properties of the $KL$ functions we can guarantee that there exists $\tau_1(\varepsilon, R, T) \geq 0$ such that $\sigma\left(b\left(R \max_{0 \leq t \leq T} \beta(t)\right), \tau_1\right) \leq \varepsilon$. Notice that if $\beta \in K^+$ then $\tau_1 \geq 0$ is independent of $T$. Then, it follows from (A.5) that we have

$$V(t) \le MAX\left\{ P\varepsilon, \Gamma(Q(\mathbf{1}\sigma(L(t_0),0))), Q\left(\mathbf{1}\zeta\left(\|u(\tau)\|_{\mathcal{U}}\big]_{t_0,t]}\right)\right)\right\} \text{ for all } t \ge t_0 + \tau_1. \text{ Since } G(s) \ge Q(\mathbf{1}\zeta(s)) \text{ for all } s \ge 0 \text{ (a}$$

consequence of (3.5)) and $Q(\mathbf{1}\varepsilon) \ge \mathbf{1}\varepsilon$, we conclude that inequality (A.2) holds for $k=1$.

Next suppose that for every $\varepsilon > 0$, $R, T \ge 0$ there exists $\tau_k(\varepsilon, R, T) \ge 0$ such that for every $(t_0, x_0, u, d) \in \mathfrak{R}_+ \times \mathcal{X} \times M_U \times M_D$ with $t_0 \in [0,T]$ and $\|x_0\|_{\mathcal{X}} \le R$ (A.2) holds for some $k \in Z_+$. Let arbitrary $\varepsilon > 0$, $R, T \ge 0$, $(t_0, x_0, u, d) \in \mathfrak{R}_+ \times \mathcal{X} \times M_U \times M_D$ with $t_0 \in [0,T]$ and $\|x_0\|_{\mathcal{X}} \le R$. Notice that the weak semigroup property implies that $\pi(t_0, x_0, u, d) \cap [t_0 + \tau_k, t_0 + \tau_k + r] \ne \varnothing$. Let $t_k \in \pi(t_0, x_0, u, d) \cap [t_0 + \tau_k, t_0 + \tau_k + r]$. Then (3.1) implies:

$$V(t) \le MAX\left\{ \mathbf{1}\sigma\left(L(t_k), t - t_k\right), \Gamma\left([V]_{[t_k,t]}\right), \mathbf{1}\gamma\left(\|u(\tau)\|_{\mathcal{U}}\big]_{[t_k,t]}\right)\right\}, \text{ for all } t \ge t_k \quad (A.7)$$

Moreover, inequality (A.2) gives:

$$[V]_{[t_k,t]} \le MAX\left\{ Q(\mathbf{1}\varepsilon), \Gamma^{(k)}(Q(\mathbf{1}\sigma(L(t_0),0))), G\left(\|u(\tau)\|_{\mathcal{U}}\big]_{[t_0,t]}\right)\right\}, \text{ for all } t \ge t_k \quad (A.8)$$

Inequality (A.6) also implies:

$$L(t_k) \le \max\left\{ \nu(t_k - t_0), c(t_0), a(R), p(Q(\mathbf{1}\sigma(L(t_0),0))), p\left(Q\left(\mathbf{1}\zeta\left(\|u(\tau)\|_{\mathcal{U}}\big]_{[t_0,t_k]}\right)\right)\right), p^u\left(\|u(\tau)\|_{\mathcal{U}}\big]_{[t_0,t_k]}\right)\right\} \quad (A.9)$$

Using (A.8) and the fact that $\Gamma(G(s)) \le G(s)$ for all $s \ge 0$ (a direct consequence of definition (3.5) and the fact that $\Gamma(Q(x)) \le Q(x)$ for all $x \in \mathfrak{R}_+^n$), we obtain:

$$\Gamma\left([V]_{[t_k,t]}\right) \le MAX\left\{ Q(\mathbf{1}\varepsilon), \Gamma^{(k+1)}(Q(\mathbf{1}\sigma(L(t_0),0))), G\left(\|u(\tau)\|_{\mathcal{U}}\big]_{[t_0,t]}\right)\right\}, \text{ for all } t \ge t_k \quad (A.10)$$

Inequality (A.10) in conjunction with inequality (A.7), the fact that $G(s) \ge Q(\mathbf{1}\zeta(s)) \ge \mathbf{1}\zeta(s)$ for all $s \ge 0$ and the fact that $t_k \le t_0 + \tau_k + r$ implies:

$$V(t) \le MAX\left\{ \mathbf{1}\sigma\left(L(t_k), t - t_0 - \tau_k - r\right), Q(\mathbf{1}\varepsilon), \Gamma^{(k+1)}(Q(\mathbf{1}\sigma(L(t_0),0))), G\left(\|u(\tau)\|_{\mathcal{U}}\big]_{[t_0,t]}\right)\right\},$$
$$\text{for all } t \ge t_0 + \tau_k + r \quad (A.11)$$

Inequality (A.9) in conjunction with the fact that $\mathbf{1}\sigma\left(p^u(s), 0\right) \le G(s)$, $\mathbf{1}\sigma\left(p(Q(\mathbf{1}\zeta(s))), 0\right) \le G(s)$ for all $s \ge 0$ and the facts that $t_k \le t_0 + \tau_k + r$, $t_0 \in [0,T]$ and $\|x_0\|_{\mathcal{X}} \le R$ implies that

$$\mathbf{1}\sigma\left(L(t_k), t - t_0 - \tau_k - r\right) \le MAX\left\{ \mathbf{1}\sigma\left(f(\varepsilon, T, R), t - t_0 - \tau_k - r\right), G\left(\|u(\tau)\|_{\mathcal{U}}\big]_{[t_0,t]}\right)\right\},$$
$$\text{for all } t \ge t_0 + \tau_k + r \quad (A.12)$$

where

$$f(\varepsilon, T, R) := \max\left\{ \max_{0 \le t \le \tau_k(\varepsilon,R,T) + r} \nu(t), \max_{0 \le t \le T} c(t), a(R), p\left(Q\left(\mathbf{1}\sigma\left(b\left(R \max_{0 \le t \le T} \beta(t)\right), 0\right)\right)\right)\right\} \quad (A.13)$$

The reader should notice that if $\beta, c \in K^+$ are bounded and $\tau_k$ is independent of $T$ then $f$ can be chosen to be independent of $T$ as well. Notice that by combining (A.11) and (A.12) we get:

$$V(t) \leq MAX\left\{ \mathbf{1}\sigma\left( f(\varepsilon,T,R), t-t_0-\tau_k -r \right), Q(\mathbf{1}\varepsilon), \Gamma^{(k+1)}\left(Q(\mathbf{1}\sigma(L(t_0),0))\right), G\left(\left[\|u(\tau)\|_{\mathcal{U}}\right]_{[t_0,t]}\right) \right\},$$

$$\text{for all } t \geq t_0 + \tau_k + r \tag{A.14}$$

Clearly, there exists $\tau(\varepsilon, R, T) \geq 0$ such that $\sigma(f(\varepsilon,T,R),\tau) \leq \varepsilon$. Define:

$$\tau_{k+1}(\varepsilon, R, T) = \tau_k(\varepsilon, R, T) + r + \tau(\varepsilon, R, T) \tag{A.15}$$

Again, the reader should notice that if $f$ and $\tau_k$ are independent of $T$ then $\tau_{k+1}$ is independent of $T$ as well. Since $Q(\mathbf{1}\varepsilon) \geq \mathbf{1}\varepsilon$ we obtain from (A.14):

$$V(t) \leq MAX\left\{ Q(\mathbf{1}\varepsilon), \Gamma^{(k+1)}\left(Q(\mathbf{1}\sigma(L(t_0),0))\right), G\left(\left[\|u(\tau)\|_{\mathcal{U}}\right]_{[t_0,t]}\right) \right\}, \text{ for all } t \geq t_0 + \tau_{k+1} \tag{A.16}$$

which shows that (A.2) holds for $k+1$.

The proof is complete.  ◁

**Proof of Theorem 4.1:** We want to show that all hypotheses of Theorem 3.1 hold with

$$L(t,x) := \max\left\{ W(t,x), \max_{i=1,\ldots,k} V_i(t,x) \right\} \tag{A.17}$$

Notice that hypothesis (H3) of Theorem 3.1 is a direct consequence of inequalities (4.2), (4.3) and definition (A.17). Moreover, hypothesis (H4) of Theorem 3.1 is a direct consequence of inequality (4.2) with $q(x) := a_1^{-1}\left( \max_{i=1,\ldots,k} x_i \right)$ for all $x \in \Re_+^n$.

Consider a solution $x(t)$ of (4.1) corresponding to arbitrary $(u,d) \in M_U \times M_D$ with arbitrary initial condition $x(t_0) = x_0 \in \Re^n$. Clearly, there exists a maximal existence time for the solution denoted by $t_{\max} \leq +\infty$. Let $V_i(t) = V_i(t,x(t))$, $i=1,\ldots,k$, $W_i(t) = W_i(t,x(t))$ absolutely continuous functions on $[t_0, t_{\max})$ and let $L(t) = L(t,x(t))$. Moreover, let $I \subset [t_0, t_{\max})$ be the zero Lebesgue measure set where $x(t)$ is not differentiable or $\dot{x}(t) \neq f(t,x(t),u(t),d(t))$. By virtue of (4.5), it follows that the following implication holds for $t \in [t_0, t_{\max}) \setminus I$ and $i=1,\ldots,k$:

$$V_i(t) \geq \max\left\{ \zeta(|u(t)|), \max_{j=1,\ldots,k} \gamma_{i,j}(V_j(t)) \right\} \Rightarrow \dot{V}_i(t) \leq -\rho_i(V_i(t)) \tag{A.18}$$

and by virtue of (4.4) we get for $t \in [t_0, t_{\max}) \setminus I$:

$$\dot{W}(t) \leq -W(t) + \lambda \max\left\{ \zeta(|u(t)|), \max_{j=1,\ldots,k} p_j(V_j(t)) \right\} \tag{A.19}$$

Lemma 3.5 in [24] in conjunction with (A.18) implies that there exists a family of continuous functions $\sigma_i$ ($i=1,\ldots,k$) of class $KL$, with $\sigma_i(s,0) = s$ for all $s \geq 0$ such that for all $t \in [t_0, t_{\max})$ and $i=1,\ldots,k$ we have:

$$V_i(t) \leq \max\left\{ \begin{array}{l} \sigma_i(V_i(t_0), t-t_0); \sup_{t_0 \leq \tau \leq t} \sigma_i\left( \max_{j=1,\ldots,k} \sup_{t_0 \leq s \leq \tau} \gamma_{i,j}(V_j(s)), t-\tau \right) \\ \sup_{t_0 \leq \tau \leq t} \sigma_i\left( \zeta\left( \sup_{t_0 \leq s \leq \tau} |u(s)| \right), t-\tau \right) \end{array} \right\} \tag{A.20}$$

Moreover, inequality (A.19) directly implies that for all $t \in [t_0, t_{\max})$ we have:

$$W(t) \leq W(t_0) + \lambda \max\left\{ \zeta\left( \sup_{t_0 \leq s \leq t} |u(s)| \right), \max_{j=1,\ldots,k} p_j\left( \sup_{t_0 \leq s \leq t} V_j(s) \right) \right\} \quad (A.21)$$

Let $\sigma(s,t) := \max_{i=1,\ldots,k} \sigma_i(s,t)$, which is a function of class $KL$ that satisfies $\sigma(s,0) = s$ for all $s \geq 0$. It follows from (A.20), (A.21) and definition (A.17) that for all $t \in [t_0, t_{max})$ and $i = 1,\ldots,k$ we get:

$$V_i(t) \leq \max\left\{ V_i(t_0), \max_{j=1,\ldots,k} \gamma_{i,j}\left( \sup_{t_0 \leq s \leq t} V_j(s) \right), \zeta\left( \sup_{t_0 \leq s \leq t} |u(s)| \right) \right\} \quad (A.22)$$

$$V_i(t) \leq \max\left\{ \sigma(L(t_0), t - t_0), \max_{j=1,\ldots,k} \gamma_{i,j}\left( \sup_{t_0 \leq s \leq t} V_j(s) \right), \zeta\left( \sup_{t_0 \leq s \leq t} |u(s)| \right) \right\} \quad (A.23)$$

$$W(t) \leq \max\left\{ \frac{1}{1-\lambda} W(t_0), \zeta\left( \sup_{t_0 \leq s \leq t} |u(s)| \right), \max_{i=1,\ldots,k} p_i\left( \sup_{t_0 \leq s \leq t} V_i(s) \right) \right\} \quad (A.24)$$

Clearly, inequalities (A.23) show that (3.1) holds with $\Gamma : \Re_+^k \to \Re_+^k$, $\Gamma(x) = (\Gamma_1(x),\ldots,\Gamma_n(x))'$ with $\Gamma_i(x) = \max_{j=1,\ldots,k} \gamma_{i,j}(x_j)$ for all $i=1,\ldots,k$ and $x \in \Re_+^n$. Furthermore, since $\gamma_{i,i}(s) \equiv 0$ for $i=1,\ldots,k$ and statement (ii) of Proposition 2.7 holds, it follows that hypothesis (H2) of Theorem 3.1 holds as well. Moreover, inequalities (A.22) and (A.24) imply that the following estimates hold for all $t \in [t_0, t_{max})$:

$$W(t) \leq \max\left\{ \begin{array}{l} \frac{1}{1-\lambda} W(t_0), \zeta\left( \sup_{t_0 \leq s \leq t} |u(s)| \right), \\ \max_{i=1,\ldots,k}\left( \max\left\{ p_i(V_i(t_0)), \max_{j=1,\ldots,k} p_i\left( \gamma_{i,j}\left( \sup_{t_0 \leq s \leq t} V_j(s) \right) \right), p_i\left( \zeta\left( \sup_{t_0 \leq s \leq t} |u(s)| \right) \right) \right\} \right) \end{array} \right\} \quad (A.25)$$

$$\max_{i=1,\ldots,k} V_i(t) \leq \max\left\{ \max_{i=1,\ldots,k} V_i(t_0), \max_{i=1,\ldots,k} \max_{j=1,\ldots,k} \gamma_{i,j}\left( \sup_{t_0 \leq s \leq t} V_j(s) \right), \zeta\left( \sup_{t_0 \leq s \leq t} |u(s)| \right) \right\} \quad (A.26)$$

Define:

$$p^u(s) := \max\left\{ \zeta(s), \max_{i=1,\ldots,k} p_i(\zeta(s)) \right\}, \text{ for all } s \geq 0 \quad (A.27)$$

$$p(x) := \max\left\{ \max_{i=1,\ldots,k} \max_{j=1,\ldots,k} \gamma_{i,j}(x_j), \max_{i=1,\ldots,k} \max_{j=1,\ldots,k} p_i(\gamma_{i,j}(x_j)) \right\}, \text{ for all } x \in \Re_+^n \quad (A.28)$$

Combining estimates (A.25), (A.26) and exploiting definitions (A.17), (A.27) and (A.28) we get for all $t \in [t_0, t_{max})$:

$$L(t) \leq \max\left\{ \frac{1}{1-\lambda} L(t_0), p\left( \sup_{t_0 \leq s \leq t} V_1(s),\ldots, \sup_{t_0 \leq s \leq t} V_k(s) \right), p^u\left( \sup_{t_0 \leq s \leq t} |u(s)| \right) \right\} \quad (A.29)$$

Inequality (3.2) is a direct consequence of (A.29), inequalities (4.2), (4.3) and Corollary 10 in [36] with $v(t) \equiv 1$, $p^u \in N_1$, $p \in N_n$ as defined by (A.27), (A.28) and appropriate $a \in N_1$ and $c \in K^+$. The reader should notice that if $\beta \in K^+$ is bounded then $c \in K^+$ is bounded as well.

Consequently, all hypotheses of Theorem 3.1 hold with $\sigma(s,t) := \max_{i=1,\ldots,k} \sigma_i(s,t)$, which is a function of class $KL$ that satisfies $\sigma(s,0) = s$ for all $s \geq 0$. The rest of proof is a consequence of Remark 3.2 in conjunction with definitions (A.27), (A.28). The proof is complete. ◁

**Proof of Theorem 4.4:** We want to show that all hypotheses of Theorem 3.1 hold with

$$L(t,x) := \max\left\{W(t,x), \max_{i=1,\ldots,k} V_i(t,x)\right\} \tag{A.30}$$

Notice that hypothesis (H3) of Theorem 3.1 is a direct consequence of inequalities (4.12), (4.13), definitions (A.30), (4.15) and Corollary 10 in [36]. Moreover, hypothesis (H4) of Theorem 3.1 is a direct consequence of inequality (4.12) with $q(x) := a_1^{-1}\left(\max_{i=1,\ldots,k} x_i\right)$ for all $x \in \Re_+^n$.

We next show that hypotheses (H1) and (H2) of Theorem 3.1 hold as well. The proof consists of two steps:

<u>Step 1:</u> We show that hypotheses (H1), (H2) of Theorem 3.1 hold for arbitrary $(t_0, u, d) \in \Re_+ \times M_U \times M_D$ and $T_r(t_0)x = x_0 \in C^1([-r,0]; \Re^n)$.

<u>Step 2:</u> We show that hypotheses (H1), (H2) of Theorem 3.1 hold for arbitrary $(t_0, u, d) \in \Re_+ \times M_U \times M_D$ and $T_r(t_0)x = x_0 \in C^0([-r,0]; \Re^n)$.

<u>Step 1:</u> Consider the solution $x(t)$ of (4.10) corresponding to arbitrary $(u,d) \in M_U \times M_D$ with arbitrary initial condition $T_r(t_0)x = x_0 \in C^1([-r,0]; \Re^n)$. Clearly, there exists a maximal existence time for the solution denoted by $t_{\max} \leq +\infty$. By virtue of Lemma A.2 in [22], we can guarantee that the functions $Q_i(t) = Q_i(t, T_{r_i}(t)x)$, $i = 1,\ldots,k$, $Q_0(t) = Q_0(t, T_{r_0}(t)x)$ are absolutely continuous functions on $[t_0, t_{\max})$. Let $V_i(t) = V_i(t, T_r(t)x) = \sup_{\theta \in [-r+r_i, 0]} Q_i(t+\theta)$, $i = 1,\ldots,k$, $W(t) = W(t, T_r(t)x) = \sup_{\theta \in [-r+r_0, 0]} Q_0(t+\theta)$ and $L(t) = L(t, T_r(t)x)$ be mappings defined on $[t_0, t_{\max})$. Moreover, let $I \subset [t_0, t_{\max})$ be the zero Lebesgue measure set where $x(t)$ is not differentiable or $\dot{x}(t) \neq f(t, T_r(t)x, u(t), d(t))$. By virtue of (4.16), it follows that the following implication holds for $t \in [t_0, t_{\max}) \setminus I$ and $i = 1,\ldots,k$:

$$Q_i(t) \geq \max\left\{\sup_{t_0 \leq s \leq t} \zeta(|u(s)|), \max_{j=1,\ldots,k} \sup_{t_0 \leq s \leq t} \gamma_{i,j}(V_j(s))\right\} \Rightarrow \dot{Q}_i(t) \leq -\rho_i(Q_i(t)) \tag{A.31}$$

and by virtue of (4.14) we have

$$\dot{Q}_0(t) \leq -Q_0(t) + \lambda \max\left\{\zeta(|u(t)|), \max_{j=1,\ldots,k} p_j(V_j(t))\right\} \tag{A.32}$$

Lemma 3.5 in [24] in conjunction with (A.31) implies that there exists a family of continuous functions $\tilde{\sigma}_i$ ($i = 1,\ldots,k$) of class $KL$, with $\tilde{\sigma}_i(s,0) = s$ for all $s \geq 0$ such that for all $t \in [t_0, t_{\max})$ and $i = 1,\ldots,k$ we have:

$$Q_i(t) \leq \max\left\{\begin{array}{l} \tilde{\sigma}_i(Q_i(t_0), t-t_0); \sup_{t_0 \leq \tau \leq t} \tilde{\sigma}_i\left(\max_{j=1,\ldots,k} \sup_{t_0 \leq s \leq \tau} \gamma_{i,j}(V_j(s)), t-\tau\right) \\ \sup_{t_0 \leq \tau \leq t} \tilde{\sigma}_i\left(\zeta\left(\sup_{t_0 \leq s \leq \tau} |u(s)|\right), t-\tau\right) \end{array}\right\} \tag{A.33}$$

Moreover, inequality (A.32) directly implies that for all $t \in [t_0, t_{\max})$ we have:

$$Q_0(t) \leq Q_0(t_0) + \lambda \max\left\{\zeta\left(\sup_{t_0 \leq s \leq t} |u(s)|\right), \max_{j=1,\ldots,k} p_j\left(\sup_{t_0 \leq s \leq t} V_j(s)\right)\right\} \tag{A.34}$$

Using the fact that $\tilde{\sigma}_i(s,0) = s$ for all $s \geq 0$, we obtain from (A.33) for all $t \in [t_0, t_{\max})$ and $i = 1,\ldots,k$:

$$Q_i(t) \leq \max\left\{\tilde{\sigma}_i(Q_i(t_0), t-t_0), \max_{j=1,\ldots,k} \gamma_{i,j}\left(\sup_{t_0 \leq s \leq t} V_j(s)\right), \zeta\left(\sup_{t_0 \leq s \leq t} |u(s)|\right)\right\} \tag{A.35}$$

Let $\sigma_i$ ($i = 1,...,k$) be functions of class $KL$, defined by $\sigma_i(s,t) = s$ for all $s \geq 0$, $t \in [0,r]$ and $\sigma_i(s,t) = \tilde{\sigma}_i(s,t-r)$ for all $s \geq 0$, $t > r$. Using the fact that $V_i(t) = V_i(t, T_r(t)x) = \sup_{\theta \in [-r+r_i, 0]} Q_i(t+\theta)$, $i = 1,...,k$ we obtain from (A.35) for all $t \in [t_0, t_{max})$ and $i = 1,...,k$:

$$V_i(t) \leq \max\left\{\sigma_i(V_i(t_0), t-t_0), \max_{j=1,...,k} \gamma_{i,j}\left(\sup_{t_0 \leq s \leq t} V_j(s)\right), \zeta\left(\sup_{t_0 \leq s \leq t} |u(s)|\right)\right\} \quad (A.36)$$

Similarly, using (A.34) and the fact that $W(t) = W(t, T_r(t)x) = \sup_{\theta \in [-r+r_0, 0]} Q_0(t+\theta)$ we obtain for all $t \in [t_0, t_{max})$:

$$W(t) \leq \max\left\{\frac{1}{1-\lambda}W(t_0), \zeta\left(\sup_{t_0 \leq s \leq t} |u(s)|\right), \max_{i=1,...,k} p_i\left(\sup_{t_0 \leq s \leq t} V_i(s)\right)\right\} \quad (A.37)$$

Let $\sigma(s,t) := \max_{i=1,...,k} \sigma_i(s,t)$, which is a function of class $KL$ that satisfies $\sigma(s,0) = s$ for all $s \geq 0$. It follows from (A.36) and definition (A.30) that for all $t \in [t_0, t_{max})$ and $i = 1,...,k$ we get:

$$V_i(t) \leq \max\left\{V_i(t_0), \max_{j=1,...,k} \gamma_{i,j}\left(\sup_{t_0 \leq s \leq t} V_j(s)\right), \zeta\left(\sup_{t_0 \leq s \leq t} |u(s)|\right)\right\} \quad (A.38)$$

$$V_i(t) \leq \max\left\{\sigma(L(t_0), t-t_0), \max_{j=1,...,k} \gamma_{i,j}\left(\sup_{t_0 \leq s \leq t} V_j(s)\right), \zeta\left(\sup_{t_0 \leq s \leq t} |u(s)|\right)\right\} \quad (A.39)$$

Clearly, inequalities (A.39) show that (3.1) holds with $\Gamma: \Re_+^k \to \Re_+^k$, $\Gamma(x) = (\Gamma_1(x),...,\Gamma_n(x))'$ with $\Gamma_i(x) = \max_{j=1,...,k} \gamma_{i,j}(x_j)$ for all $i = 1,...,k$ and $x \in \Re_+^n$. Moreover, since statement (ii) of Proposition 2.7 holds, it follows that hypothesis (H2) of Theorem 3.1 holds as well. Finally, inequalities (A.37) and (A.38) imply that the following estimates hold for all $t \in [t_0, t_{max})$:

$$W(t) \leq \max\left\{\begin{array}{l} \dfrac{1}{1-\lambda}W(t_0), \zeta\left(\sup_{t_0 \leq s \leq t} |u(s)|\right), \\ \max_{i=1,...,k}\left(\max\left\{p_i(V_i(t_0)), \max_{j=1,...,k} p_i\left(\gamma_{i,j}\left(\sup_{t_0 \leq s \leq t} V_j(s)\right)\right), p_i\left(\zeta\left(\sup_{t_0 \leq s \leq t} |u(s)|\right)\right)\right\}\right) \end{array}\right\} \quad (A.40)$$

$$\max_{i=1,...,k} V_i(t) \leq \max\left\{\max_{i=1,...,k} V_i(t_0), \max_{i=1,...,k} \max_{j=1,...,k} \gamma_{i,j}\left(\sup_{t_0 \leq s \leq t} V_j(s)\right), \zeta\left(\sup_{t_0 \leq s \leq t} |u(s)|\right)\right\} \quad (A.41)$$

Define:

$$p^u(s) := \max\left\{\zeta(s), \max_{i=1,...,k} p_i(\zeta(s))\right\}, \text{ for all } s \geq 0 \quad (A.42)$$

$$p(x) := \max\left\{\max_{i=1,...,k} \max_{j=1,...,k} \gamma_{i,j}(x_j), \max_{i=1,...,k} \max_{j=1,...,k} p_i(\gamma_{i,j}(x_j))\right\}, \text{ for all } x \in \Re_+^n \quad (A.43)$$

Combining estimates (A.40), (A.41) and exploiting definitions (A.30), (A.42) and (A.43) we get for all $t \in [t_0, t_{max})$:

$$L(t) \leq \max\left\{\frac{1}{1-\lambda}L(t_0), p\left(\sup_{t_0 \leq s \leq t} V_1(s),..., \sup_{t_0 \leq s \leq t} V_k(s)\right), p^u\left(\sup_{t_0 \leq s \leq t} |u(s)|\right)\right\} \quad (A.44)$$

Inequality (3.2) is a direct consequence of (A.44), inequalities (4.12), (4.13) and Corollary 10 in [36] with $v(t) \equiv 1$, $p^u \in N_1$, $p \in N_n$ as defined by (A.42), (A.43) and appropriate $a \in N_1$ and $c \in K^+$. The reader should notice that if $\beta \in K^+$ is bounded then $c \in K^+$ is bounded as well.

**Step 2:** Let $(t_0, x_0, u, d) \in \Re_+ \times C^1([-r,0]; \Re^n) \times M_U \times M_D$. Inequalities (A.39) imply for the solution $x(t)$ of (4.10) corresponding to $(u,d) \in M_U \times M_D$ with initial condition $T_r(t_0)x = x_0 \in C^1([-r,0]; \Re^n)$ and for all $t \in [t_0, t_{\max})$:

$$V(t) \leq MAX \left\{ Q(1\sigma(L(t_0),0)), Q\left(1\zeta\left(\|u(\tau)\|_U\right]_{[t_0,t]}\right)\right) \right\} \quad (A.45)$$

where $Q(x) = MAX\{x, \Gamma(x), \Gamma^{(2)}(x), ..., \Gamma^{(n-1)}(x)\}$. Using (4.12), (4.13), (A.30), (A.44) and (A.45) we obtain functions $\rho \in K^+$, $a \in K_\infty$ such that the solution $x(t)$ of (4.10) corresponding to $(u,d) \in M_U \times M_D$ with initial condition $T_r(t_0)x = x_0 \in C^1([-r,0]; \Re^n)$ and for all $t \in [t_0, t_{\max})$:

$$\|T_r(t)x\|_r \leq a\left(\rho(t) + \|x_0\|_r + \sup_{t_0 \leq s \leq t} |u(s)|\right) \quad (A.46)$$

Lemma 2.6 in [23] and (A.46) imply that system (4.10) is RFC from the input $u \in M_U$.

We next claim that inequalities (3.1) and (3.2) hold for all $(t_0, x_0, u, d) \in \Re_+ \times C^0([-r,0]; \Re^n) \times M_U \times M_D$ and $t \geq t_0$. The proof will be made by contradiction. Suppose on the contrary that there exists $(t_0, x_0, u, d) \in \Re_+ \times C^0([-r,0]; \Re^n) \times M_U \times M_D$ and $t_1 > t_0$ such that the solution $x(t)$ of (4.10) with initial condition $T_r(t_0)x = x_0$ corresponding to input $(u,d) \in M_U \times M_D$ satisfies $\beta(t_1, t_0, x_0, d, u) > 0$, where

$$\beta(t, t_0, x_0, d, u) := \max \begin{cases} L(t) - \max\left\{\nu(t-t_0), c(t_0), a(\|x_0\|_r), p([V]_{[t_0,t]}), p^u\left(\|u(\tau)\|_U\right]_{[t_0,t]}\right)\right\} \\ \max_{i=1,...,k}\left\{V_i(t) - \max\left\{\sigma(L(t_0), t-t_0), \Gamma_i([V]_{[t_0,t]}), \zeta\left(\|u(\tau)\|_U\right]_{[t_0,t]}\right)\right\}\right\} \end{cases}$$

Using continuity of the mappings $x \to L(t,x)$, $x \to V_i(t,x)$ ($i=1,...,k$) and continuity of the solution of (4.10) with respect to the initial condition we can guarantee that the mapping $x_0 \to \beta(t_1, t_0, x_0, d, u)$ is continuous. Using density of $C^1([-r,0]; \Re^n)$ in $C^0([-r,0]; \Re^n)$, continuity of the mapping $x_0 \to \beta(t_1, t_0, x_0, d, u)$, we conclude that there exists $\tilde{x}_0 \in C^1([-r,0]; \Re^n)$ such that:

$$|\beta(t_1, t_0, x_0, d, u) - \beta(t_1, t_0, \tilde{x}_0, d, u)| \leq \frac{1}{2}\beta(t_1, t_0, x_0, d, u)$$

Thus we obtain a contradiction.

Consequently, all hypotheses of Theorem 3.1 hold with $\sigma(s,t) := \max_{i=1,...,k} \sigma_i(s,t)$, which is a function of class $KL$ that satisfies $\sigma(s,0) = s$ for all $s \geq 0$. The rest of proof is a consequence of Remark 3.2 in conjunction with definitions (A.42), (A.43). The proof is complete. ◁

**Proof of Theorem 4.8:** We want to show that all hypotheses of Theorem 3.1 hold with

$$L(t,x) := \max\left\{W(t,x), \max_{i=1,...,k} V_i(x)\right\} \quad (A.47)$$

where

$$W(t,x) := \exp(-(\mu+1)t)Q(x) \quad (A.48)$$

Notice that definition (A.48) in conjunction with inequalities (4.23) imply the following inequality for all $(t,x) \in \Re_+ \times \Re^n$:

$$\exp(-(\mu+1)t)a_3(|x|) - g(V_1(x),...,V_k(x)) - \kappa \leq W(t,x) \leq a_4(|x|)$$

Using Corollary 10 in [36], we can find functions $\tilde{a} \in K_\infty$, $\eta \in K^+$ such that $a_3^{-1}(s\exp((\mu+1)t)) \leq \frac{1}{\eta(t)}\tilde{a}(s)$ for all $t, s \geq 0$. Consequently, we obtain $\exp(-(\mu+1)t)a_3(s) \geq \tilde{a}^{-1}(\eta(t)s)$ for all $t, s \geq 0$. Notice that hypothesis (H3) of Theorem 3.1 with $\beta(t) \equiv 1$ is a direct consequence of previous inequalities, (4.22) and definitions (A.47), (A.48). Moreover, hypothesis (H4) of Theorem 3.1 is a direct consequence of inequality (4.22) with $q(x) := a_1^{-1}\left(\max_{i=1,\ldots,k} x_i\right)$ for all $x \in \Re_+^n$.

Consider the solution $x(t)$ of (4.21) under hypotheses (R1-4) corresponding to arbitrary $(u, d, \tilde{d}) \in M_U \times M_D \times M_{\Re_+}$ with arbitrary initial condition $x(t_0) = x_0 \in \Re^n$. Notice that since system (4.21) is autonomous (see [18]) it suffices to consider the case $t_0 = 0$. By virtue of Proposition 2.5 in [18], there exists a maximal existence time for the solution denoted by $t_{\max} \leq +\infty$. Let $V_i(t) = V_i(x(t))$, $i = 1,\ldots,k$, $W(t) = W(t, x(t))$, $L(t) = L(t, x(t))$ absolutely continuous functions on $[0, t_{\max})$. Moreover, let $\pi := \{\tau_0, \tau_1, \ldots\}$ be the set of sampling times (which may be finite if $t_{\max} < +\infty$) and $p(t) := \max\{\tau \in \pi : \tau \leq t\}$, $q(t) := \min\{\tau \in \pi : \tau \geq t\}$. Let $I \subset [0, t_{\max})$ be the zero Lebesgue measure set where $x(t)$ is not differentiable or where $\dot{x}(t) \neq f(x(t), x(\tau_i), d(t), u(t), u(\tau_i))$. Clearly, we have $x(t) = \phi(t - p(t), x(p(t)); P_t d, P_t u)$ for all $t \in [0, t_{\max})$, where $(P_t u)(s) = u(p(t) + s)$, $(P_t d)(s) = d(p(t) + s)$, $s \geq 0$. Next we show that the following implication holds for $t \in [0, t_{\max}) \setminus I$ and $i = 1,\ldots,k$:

$$V_i(t) \geq \max\left\{\zeta\left(\sup_{p(t) \leq s \leq t} |u(s)|\right), \max_{j=1,\ldots,k} \gamma_{i,j}\left(\sup_{p(t) \leq s \leq t} V_j(s)\right)\right\} \Rightarrow \dot{V}_i(t) \leq -\rho_i(V_i(t)) \tag{A.49}$$

In order to prove implication (A.49) let $t \in [0, t_{\max}) \setminus I$, $i = 1,\ldots,k$, $\tau = p(t)$ and suppose that $V_i(t) \geq \max\left\{\zeta\left(\sup_{p(t) \leq s \leq t} |u(s)|\right), \max_{j=1,\ldots,k} \gamma_{i,j}\left(\sup_{p(t) \leq s \leq t} V_j(s)\right)\right\}$. By virtue of the semigroup property for the previous inequality implies that $\zeta(|u(\tau + s)|) = \zeta(|(P_t u)(s)|) \leq V_i(x(t))$, $\gamma_{i,j}(V_j(\phi(s, x(\tau); P_t d, P_t u))) \leq V_i(x(t))$ for all $s \in [0, t - \tau]$ and $j = 1,\ldots,k$. In this case, by virtue of definition (4.26) and the fact that $t - \tau \leq h(x(\tau), u(\tau))$, it follows that $x(\tau) \in A_i(h(x(\tau), u(\tau)), x(t))$. Since $\dot{x}(t) = f(x(t), x(\tau), d(t), u(t), u(\tau))$, we conclude from (4.25) that $\dot{V}_i(t) \leq -\rho_i(V_i(t))$.

Lemma 3.5 in [24] implies that there exists a family of continuous function $\sigma_i$ of class $KL$ ($i = 1,\ldots,k$), with $\sigma_i(s, 0) = s$ for all $s \geq 0$ such that for all $t \in [0, t_{\max})$ and $i = 1,\ldots,k$ we have:

$$V_i(t) \leq \max\left\{\sigma_i(V_i(0), t); \max_{j=1,\ldots,k} \sup_{0 \leq \tau \leq t} \sigma_i\left(\gamma_{i,j}\left(\sup_{p(\tau) \leq s \leq \tau} V_j(s)\right), t - \tau\right); \sup_{0 \leq \tau \leq t} \sigma_i\left(\zeta\left(\sup_{p(\tau) \leq s \leq \tau} |u(s)|\right), t - \tau\right)\right\} \tag{A.50}$$

Let $\sigma(s, t) := \max_{i=1,\ldots,k} \sigma_i(s, t)$, which is a function of class $KL$ that satisfies $\sigma(s, 0) = s$ for all $s \geq 0$. Inequalities (3.1) with $\Gamma : \Re_+^k \to \Re_+^k$, $\Gamma(x) = (\Gamma_1(x), \ldots, \Gamma_n(x))'$ with $\Gamma_i(x) = \max_{j=1,\ldots,k} \gamma_{i,j}(x_j)$ for all $i = 1,\ldots,k$ and $x \in \Re_+^n$ are direct consequences of the previous definition, estimates (A.50), definition (A.47) and the fact that $\sigma_i(s, 0) = s$ for all $s \geq 0$ and $i = 1,\ldots,k$. Moreover, since statement (ii) of Proposition 2.7 holds, it follows that hypothesis (H2) of Theorem 3.1 holds as well.

Exploiting (4.24) and definition (A.48) we get for $t \in [0, t_{\max}) \setminus I$:

$$\dot{W}(t) \leq -W(t) + \lambda \max\left\{\zeta\left(\sup_{p(t) \leq s \leq t} |u(s)|\right), \max_{j=1,\ldots,k} p_j\left(\sup_{p(t) \leq s \leq t} V_j(s)\right)\right\} \tag{A.51}$$

Inequality (A.51) directly implies that for all $t \in [0, t_{\max})$ we have:

$$W(t) \leq \max\left\{\frac{1}{1-\lambda}W(0), \zeta\left(\sup_{0\leq s\leq t}|u(s)|\right), \max_{i=1,\ldots,k} p_i\left(\sup_{0\leq s\leq t} V_i(s)\right)\right\} \quad (A.52)$$

Moreover, inequalities (A.50) and (A.52) imply that the following estimates hold for all $t \in [0, t_{\max})$:

$$W(t) \leq \max\left\{\begin{array}{l} \frac{1}{1-\lambda}W(0), \zeta\left(\sup_{0\leq s\leq t}|u(s)|\right), \\ \max_{i=1,\ldots,k}\left(\max\left\{p_i(V_i(0)), \max_{j=1,\ldots,k} p_i\left(\gamma_{i,j}\left(\sup_{0\leq s\leq t} V_j(s)\right)\right), p_i\left(\zeta\left(\sup_{0\leq s\leq t}|u(s)|\right)\right)\right\}\right) \end{array}\right\} \quad (A.53)$$

$$\max_{i=1,\ldots,k} V_i(t) \leq \max\left\{\max_{i=1,\ldots,k} V_i(0), \max_{i=1,\ldots,k}\max_{j=1,\ldots,k} \gamma_{i,j}\left(\sup_{0\leq s\leq t} V_j(s)\right), \zeta\left(\sup_{0\leq s\leq t}|u(s)|\right)\right\} \quad (A.54)$$

Define:

$$p^u(s) := \max\left\{\zeta(s), \max_{i=1,\ldots,k} p_i(\zeta(s))\right\}, \text{ for all } s \geq 0 \quad (A.55)$$

$$p(x) := \max\left\{\max_{i=1,\ldots,k}\max_{j=1,\ldots,k} \gamma_{i,j}(x_j), \max_{i=1,\ldots,k}\max_{j=1,\ldots,k} p_i(\gamma_{i,j}(x_j))\right\}, \text{ for all } x \in \Re_+^n \quad (A.56)$$

Combining estimates (A.53), (A.54) and exploiting definitions (A.55), (A.56) and (A.47) we get for all $t \in [0, t_{\max})$:

$$L(t) \leq \max\left\{\frac{1}{1-\lambda}L(0), p\left(\sup_{0\leq s\leq t} V_1(s),\ldots,\sup_{0\leq s\leq t} V_k(s)\right), p^u\left(\sup_{0\leq s\leq t}|u(s)|\right)\right\} \quad (A.57)$$

Inequality (3.2) is a direct consequence of (A.57), inequalities (4.22), (4.23) with $\nu(t) = c(t) \equiv 1$, $p^u \in N_1$, $p \in N_n$ as defined by (A.55), (A.56) and appropriate $a \in N_1$.

Consequently, all hypotheses of Theorem 3.1 hold with $\sigma(s,t) := \max_{i=1,\ldots,k} \sigma_i(s,t)$, which is a function of class $KL$ that satisfies $\sigma(s,0) = s$ for all $s \geq 0$. The rest of proof is a consequence of Remark 3.2 in conjunction with definitions (A.55), (A.56). The proof is complete. ◁